\documentclass[reqno,11pt]{amsart}
\usepackage{xcolor}
\usepackage[top=2.0cm,bottom=2.0cm,left=3cm,right=3cm]{geometry}
\usepackage{amsthm,amsmath,amssymb,dsfont}
\usepackage{mathrsfs,amsfonts,functan,extarrows,mathtools}

\usepackage{xcolor}
\usepackage{cite}
\usepackage{indentfirst, latexsym, amssymb, enumerate,amsmath,graphicx}
\usepackage{float}
\usepackage{relsize}
\usepackage{marginnote}
\usepackage{stmaryrd}
\usepackage{esint}
\usepackage{graphicx}
\usepackage{bm}
\usepackage{caption}
\usepackage{subfigure}
\usepackage{cite}
\usepackage{color}
\usepackage{graphicx}
\usepackage{subfigure}
\usepackage{float}
\usepackage{paralist}
\usepackage{indentfirst}
\usepackage{cite}
\allowdisplaybreaks 
\theoremstyle{plain}
\newtheorem{thm}{Theorem}[section]

\newtheorem{lem}[thm]{Lemma}
\newtheorem{prop}[thm]{Proposition}
\newtheorem{rem}{Remark}[section]

\newtheorem{defn}{Definition}[section]

\numberwithin{equation}{section}

\usepackage{appendix}
\usepackage{xcolor}

\newcommand{\eps}{{\varepsilon}}

\begin{document}

\title[The compressible pressureless Navier-Stokes  system]
{Uniform stability and optimal time decay rates  of the compressible pressureless Navier-Stokes system in  the critical regularity framework}

\author[F. Li]{Fucai Li}
\address{School of Mathematics, Nanjing University, Nanjing
 210093, P. R. China}
\email{fli@nju.edu.cn}

\author[J. Ni]{Jinkai Ni} 
\address{School of Mathematics, Nanjing University, Nanjing
 210093, P. R. China}
\email{jinkaini123@gmail.com}

\author[Z. Zhang]{Zhipeng Zhang}  
\address{School of Mathematical Sciences, Ocean University of China, Qingdao
 266100, P. R. China}  
\email{zhangzp@ouc.edu.cn}

\begin{abstract}
This paper investigates the Cauchy problem for the compressible pressureless Navier-Stokes system in $\mathbb{R}^d$ with $d \geq 2$. 
Unlike the standard isentropic compressible Navier-Stokes system, the density in the pressureless model lacks a dissipative mechanism, leading to significant coupling effects from nonlinear terms in the momentum equations. 
We first prove the global well-posedness and uniform stability of strong solutions to the compressible pressureless Navier-Stokes system in the critical Besov space $\dot{B}_{2,1}^{\frac{d}{2}} \times \dot{B}_{2,1}^{\frac{d}{2}-1}$. Then, under the additional assumption that the  low-frequency component of the  initial density belongs to $\dot{B}_{2,\infty}^{\sigma_0+1}$ and  that the initial velocity  is sufficiently small in $\dot{B}_{2,\infty}^{\sigma_0}$ with $\sigma_0 \in (-\frac{d}{2}, \frac{d}{2}-1]$, we overcome the challenge of derivative loss caused by nonlinearity and establish optimal decay estimates for $u$ in $\dot{B}_{2,1}^{\sigma}$ with $\sigma \in (\sigma_0, \frac{d}{2}+1]$. 
In particular, it is shown that the density remains uniformly bounded in time which reveals a new asymptotic behavior in contrast to the isentropic compressible Navier-Stokes system where the density  exhibits a dissipative structure and decays over time.
\end{abstract}



\maketitle

\setcounter{equation}{0}
 \indent \allowdisplaybreaks

\section{Introduction}\label{Sec:intro-resul}
\subsection{ The system and some related results}
In this paper, we study the following compressible pressureless Navier-Stokes system: 
\begin{equation}\label{A1}
\left\{
\begin{aligned}
&\partial_t \rho+{\rm div}\,(\rho u)=0                    ,\\
& \partial_t( \rho u)+{\rm div}\,(\rho u \otimes u)-\mu\Delta u
 -(\mu+\nu)\nabla{\rm div}\,u =  0, 
\end{aligned}\right.
\end{equation}
in $\mathbb{R}^d$ with $d\geq 2$, with the initial data
\begin{align}\label{A2}
( \rho(t,x),u(t,x) )|_{t=0}=& ( \rho_0(x),u_0(x) ),
\end{align} 
where the unknown functions $\rho(t,x)\geq 0$ and $u(t,x)=(u_1(t,x), \dots, u_d(t,x))\in \mathbb{R}^d$ stand for the density and velocity of the fluid, respectively. 
The viscosity coefficients $\mu$ and $\nu$ satisfy $\mu > 0$ and  $2\mu +\nu>0$. 
The far-field conditions are given by
\begin{align}\label{A3}
( \rho_0(x), u_0(x) ) \rightarrow (\bar{\rho}, 0), \quad \text{as} \quad |x| \rightarrow \infty,
\end{align}
where $\bar{\rho} > 0$ denotes the constant background density.
The system \eqref{A1} is commonly used to describe various astrophysical phenomena \cite{DA}. 
Formally, it can be derived from Boltzmann-type equations that model interacting agents. This is achieved by first performing a mean-field limit to obtain the Cucker-Smale system, and then determining the evolution of momentum and mass at every point by averaging over the kinetic velocity \cite{CFTV}.
Alternatively, we can also derive the system \eqref{A1}   from the isentropic compressible Navier-Stokes system by taking the high Mach number limit.
However, only a few of results are available for the high Mach number limit, 
see \cite{Haspot-2016-JMFM, Haspot-2016-DCDS,L}. 
One major challenge arises from the lack of compactness in the density, as its bound is no longer preserved for the constant viscous coefficients case.

Given the close relationship between the compressible Navier-Stokes system and the system \eqref{A1}, we first recall some progress made on the global well-posedness of the former. 
The global well-posedness of classical solutions with the  initial data close to a non-vacuum equilibrium in Sobolev space $H^s$ was first obtained by Matsumura and Nishida \cite{MN-JMKU-1980}. Such a theory was later generalized to weak solutions by Hoff \cite{Hoff-1995-JDE}.
It should be noted that there have been extensive studies on the solvability problem in $\mathbb R^d$ with $d\geq 2$ in the so-called ``critical regularity" framework. The central idea, originating from Fujita and Kato's work \cite{FK-1964-ARMA} on the incompressible homogeneous Navier-Stokes equations, is that the “optimal” function spaces for the well-posedness of the compressible Navier-Stokes system 
must be invariant under the following scaling transformations for any $\lambda>0$:  
\begin{align*}
\rho_{\lambda}(t,x) = \rho(\lambda^2 t, \lambda x), \quad u_{\lambda}(t,x) = \lambda u(\lambda^2 t, \lambda x),
\end{align*}
along with the corresponding scaling of the initial data:
\begin{align*}
\rho_0(x) \leadsto \rho_0(\lambda x), \quad u_0(x) \leadsto \lambda u_0(\lambda x).
\end{align*}
The global well-posedness of the Cauchy problem was first investigated by Danchin \cite{Danchin-IM-00} in $L^2$-type Besov spaces, and subsequently extended by Charve and Danchin \cite{CD-2010-ARMA}, Chen et al. \cite{CMZ-2010-CPAM} and Haspot \cite{Haspot-2011-ARMA} to general $L^p$-type Besov spaces. Furthermore, significant progress has been made on optimal time-decay estimates for solutions to the Cauchy problem in critical $L^2$-type and $L^p$-type Besov spaces, see   \cite{Danchin-CPDE-2007,DX-ARMA-2017,Xj-CMP-2019,X-X-2021-JDE,Haspot-2011-ARMA}.
For the existence of solutions for large data, the major breakthrough is due to Lions \cite{Lions-1998} and Feireisl et al. \cite{Feireisl-2001-JMFM}, where they obtained the global existence of weak solutions when the adiabatic exponent is suitably large.  
Recently, Huang, Li and Xin \cite{HLX-2012-CPAM} established the global existence and uniqueness of classical solutions to the Cauchy problem for the isentropic compressible Navier-Stokes equations in three-dimensional space with smooth initial data which are of small energy but possibly large oscillations; in particular, the initial density is allowed to vanish and even has compact support.

However, as Danchin pointed out in \cite{Dr-2024-PMP}, unlike the isentropic compressible Navier-Stokes equations, the system \eqref{A1} lacks a pressure term, making it difficult to control the density or use the properties of the standard transport equation. In particular, the compactness methods of   Lions \cite{Lions-1998} and  Feireisl et al. \cite{Feireisl-2001-JMFM} for viscous isentropic compressible fluids with the pressure law $P(\rho)=a\rho^\gamma$ (where $a>0$, $\gamma>3/2$) is invalid here, as there is no {\it viscous effective flux} analogue to establish the compactness needed to pass to the limit in \eqref{A1} from approximate solutions.
As a result, studies on the global well-posedness of the system \eqref{A1} are limited.
Danchin \cite{Dr-2024-PMP} established a global existence and uniqueness result for the problem \eqref{A1}--\eqref{A2} with large density variations in $\mathbb{R}^2$. Guo et al. \cite{GTZ-JDE-2025} proved the global existence of the problem \eqref{A1}--\eqref{A2} in $\mathbb{R}^3$, and achieved the optimal time decay rate for the second-order derivative of $u$, under the assumptions that the initial datum for $\rho_0$ has $H^3$ regularity, while   for $u_0$ satisfies $H^4$ regularity, along with a smallness condition in $L^1$. Recently, by employing the method developed in \cite{Zp-ADV-2020,HSWZ-2025-IMRN}, Wang et al. \cite{WWX-2025-arXiv} investigated the $\dot B_{2,1}^{\frac{1}{2}}$-type Fujita-Kato solution to the problem \eqref{A1}--\eqref{A2} in $\mathbb{R}^3$ under the assumption that the initial density $\rho_0$ is merely bounded. Moreover, they obtained the decay estimate for the upper bound of $u$ in $\dot B_{2,1}^{\frac{3}{2}}$, namely, the decay rate of the $1.5$-th order derivative of $u$.

However, the uniform stability of the problem \eqref{A1}--\eqref{A2} under the lower critical regularity assumption and the optimal decay rate of $u$ have not yet been fully established. A more significant challenge stems from the lack of dissipative structure in $\rho$. The nonlinear terms in the momentum equations involve factors of $\frac{1}{\rho}$, such as $\frac{\Delta u}{\rho}$ and $\frac{\nabla{\rm div}\, u}{\rho}$, which result in a loss of derivatives in the solution $u$ and pose substantial difficulties in deriving decay estimates. 
In \cite{WWX-2025-arXiv}, the boundedness assumption on the initial density $\rho_0$  helps overcome these issues. By assuming the initial density belongs to the critical regularity space $\dot B_{2,1}^{\frac{d}{2}}$, we can more precisely characterize the adverse effects induced by these coupled nonlinear terms, thereby achieving a deeper understanding of the solution behavior for the problem \eqref{A1}--\eqref{A2}. The main objective of this work is to address these unresolved problems.
In the present paper, we first establish the existence of global strong solutions under weaker critical regularity assumptions in the homogeneous Besov space $\dot B_{2,1}^{\frac{d}{2}} \times \dot B_{2,1}^{\frac{d}{2}-1}$. After that, we prove uniform stability and derive optimal decay rates for a broader class of derivatives of $u$ up to order $2.5$ (for $d=3$), thereby improving the decay results previously established in \cite{GTZ-JDE-2025,WWX-2025-arXiv} in $\mathbb R^3$.  

\subsection{Notations and basic facts}
Before stating our main results, we introduce some   notations and the definitions of Besov spaces used throughout this paper.
The symbol $C$ denotes a generic positive constant that may vary from line to line. The notation $A \lesssim B$ (resp. $A \gtrsim B$) means $A \leq CB$ (resp. $A \geq CB$) for some constant $C > 0$. For any Banach space $X$ and functions $f, g \in X$, we define $\|(f,g)\|_{X} := \|f\|_{X} + \|g\|_{X}$. For any $T > 0$ and $1 \leq \varrho \leq \infty$, $L^\varrho(0,T;X)$ denotes the space of measurable functions $f: [0,T] \to X$ such that the mapping $t \mapsto \|f(t)\|_{X}$ belongs to $L^\varrho(0,T)$, equipped with the norm $\|\cdot\|_{L^\varrho(0,T;X)} = \|\cdot\|_{L_T^\varrho(X)}$. $\mathcal{C}([0,T];X)$ denotes the space of continuous functions $f:[0,T]\rightarrow X$. Let $\mathcal{F}(f):=\widehat{f}$ and $\mathcal{F}^{-1}(f):=\check{f}$ denote the Fourier transform of $f$ and its inverse, respectively.

We next recall the Littlewood-Paley decomposition and the definitions of Besov spaces; for further details,  interested reader can refer to \cite[Chapters 2–3]{BCD-Book-2011}. Let $\chi(\xi)$ be a smooth, radial, non-increasing function supported in $B(0,\frac{4}{3})$ such that $\chi(\xi) = 1$ on $B(0,\frac{3}{4})$. Then the function $\phi(\xi) := \chi(\frac{\xi}{2}) - \chi(\xi)$ satisfies
\begin{align*}
\sum_{k\in\mathbb{Z}} \phi(2^{-k}\cdot) = 1,
\quad \text{and} \quad
\mathrm{supp}\,\phi \subset \Big\{ \xi \in \mathbb{R}^d \, \Big|\, \frac{3}{4} \leq |\xi| \leq \frac{8}{3} \Big\}.
\end{align*}
For each $k \in \mathbb{Z}$, the homogeneous dyadic block $\dot{\Delta}_k$ is defined by
\begin{align*}
\dot{\Delta}_k f := \mathcal{F}^{-1}\big( \phi(2^{-k}\cdot)\mathcal{F}(f) \big) = 2^{kd} h(2^k\cdot) \ast f,
\quad \text{with} \quad h := \mathcal{F}^{-1}\phi.
\end{align*}
Let $\mathcal{P}$ denote the class of all polynomials  on $\mathbb R^d$ and
$\mathcal{S}^\prime_h=\mathcal{S}^\prime/\mathcal{P}$ represent the tempered distributions on $\mathbb R^d$ modulo polynomials.  Then for any $f\in \mathcal{S}^\prime_h$, one has 
\begin{align*}
f=\sum_{k\in\mathbb Z}\dot\Delta_k f\quad\text{ for any}\quad f\in \mathcal{S}^\prime_h,\quad \dot\Delta_k\dot\Delta_j=0\quad {\rm if}\quad |k-j|\geq2.  \end{align*}

Thanks to those dyadic blocks, Besov spaces are defined as follows.
\begin{defn}
For $s\in \mathbb R$, and $1\leq p,r\leq\infty$, the homogeneous Besov spaces $\dot B_{p,r}^s$ are defined by  
\begin{align*}
\dot B_{p,r}^s:=\big\{ f\in \mathcal{S}^\prime_h \,\big{|}\, \|f\|_{\dot{B}_{p,r}^s}:=\big\|\{2^{ks}\|\dot\Delta_k u\|_{L^p}\}_{k\in\mathbb Z}\big\|_{l^r}<\infty\big\}. 
\end{align*}
\end{defn}
Furthermore, we recall   a class of mixed space-time Besov spaces which  was originally introduced by  Chemin and  Lerner \cite{C-N-1995} (see also the special case of Sobolev spaces in \cite{Chemin-1999}).

\begin{defn}
For $T>0$, $s\in \mathbb R$, $1\leq \varrho,p,r\leq \infty $,  the homogeneous Chemin-Lerner space  $\widetilde L^\varrho(0,T;\dot B_{p,r}^s)$ is defined by
\begin{align*}
\widetilde L^\varrho(0,T;\dot B_{p,r}^s):=\big\{f\in L^\varrho(0,T;\mathcal{S}^\prime_h)\,\big{|}\, \|f\|_{\widetilde L^\varrho(0,T;\dot B_{p,r}^s)}:= \big\|\{2^{ks}\|\dot\Delta_k f\|_{L^\varrho_T(L^p)}\}_{k\in\mathbb Z}\big\|_{l^r}<\infty             \big\} .   
\end{align*}
\end{defn}
By applying Minkowski's inequality, one gets the following facts: 
\begin{rem}
It holds that
\begin{align*}
\|f\|_{\widetilde L^\varrho_T(\dot B_{p,r}^s)}\leq \|f\|_{  L^\varrho_T(\dot B_{p,r}^s)}\quad {\rm if} \quad r\geq \varrho; \quad \|f\|_{\widetilde L^\varrho_T(\dot B_{p,r}^s)}\geq \|f\|_{  L^\varrho_T(\dot B_{p,r}^s)}\quad {\rm if} \quad r\leq \varrho.        
\end{align*}
Here, $\|\cdot\|_{L_T^\varrho(\dot B_{p,r}^s)}$ is the usual Lebesgue-Besov norm.
\end{rem}
Restricting the Besov norms to the low- or high-frequency components of distributions plays a crucial role in our approach. We frequently employ the following notations for any $s\in \mathbb R$ and $1\leq p,r\leq \infty$:
\begin{gather*} 
    \|f\|^\ell_{\dot B_{p,r}^s}:=\big\|\{2^{ks}\|\dot\Delta_k f\|_{L^p}\}_{k\leq 0} \big\|_{l^r}, \quad\|f\|^\ell_{\widetilde L^\varrho_T(\dot B_{p,r}^s)}:=\big\|\{2^{ks}\|\dot\Delta_k f\|_{L^\varrho_T(L^p)}\}_{k\leq 0} \big\|_{l^r},      \\
    \|f\|^h_{\dot B_{p,r}^s}:=\big\|\{2^{ks}\|\dot\Delta_k f\|_{L^p}\}_{k\geq -1} \big\|_{l^r},\quad  \|f\|^h_{\widetilde L^\varrho_T(\dot B_{p,r}^s)}:=\big\|\{2^{ks}\|\dot\Delta_k f\|_{L^\varrho_T(L^p)}\}_{k\geq -1} \big\|_{l^r}.            
\end{gather*}
Define
\begin{align*}
 f ^\ell:= \sum_{k\leq -1}\dot\Delta_k f,\quad  f ^h:= \sum_{k\geq 0}\dot\Delta_k f.  
\end{align*}
It is evident for any $s^\prime > 0$ that
\begin{gather*}
    \|f^\ell\|_{\dot B_{p,r}^s}\lesssim \|f\|_{\dot B_{p,r}^s}^\ell\lesssim \|f\|_{\dot B_{p,r}^{s-s^\prime}}^\ell, \quad\|f^h\|_{\dot B_{p,r}^s}\lesssim \|f\|_{\dot B_{p,r}^s}^h\lesssim \|f\|_{\dot B_{p,r}^{s+s^\prime}}^h,     \\
   \|f^\ell\|_{\widetilde L^\varrho_T(\dot B_{p,r}^s)}\lesssim \|f\|_{\widetilde L^\varrho_T(\dot B_{p,r}^s)}^\ell\lesssim \|f\|_{\widetilde L^\varrho_T(\dot B_{p,r}^{s-s^\prime})}^\ell,  \quad \|f^h\|_{\widetilde L^\varrho_T(\dot B_{p,r}^s)}\lesssim \|f\|_{\widetilde L^\varrho_T(\dot B_{p,r}^s)}^h\lesssim \|f\|_{\widetilde L^\varrho_T(\dot B_{p,r}^{s+s^\prime})}^h.          
\end{gather*}

\subsection{Main results} With the above preparations in hand, we now present our main results.
Without loss of generality, we set  $\bar{\rho}=1$.
By defining the fluctuation variables $ a := \rho - 1 $ and $ a_0 := \rho_0 - 1 $, the Cauchy problem \eqref{A1}--\eqref{A2} is reformulated as 
\begin{equation}\label{A4}\left\{
\begin{aligned}
&\partial_t a+{\rm div}\, u=-{\rm div}\,(au)                   ,\\
& \partial_t u-\mu\Delta u-(\mu+\nu)\nabla {\rm div}\, u = -u\cdot\nabla u+ \mu f(a)\Delta u+(\mu+\nu)f(a)\nabla{\rm div}\, u, \\
& (a,u)(x,0)=(a_0,u_0)(x)\rightarrow (0,0),\quad |x|\rightarrow \infty,
\end{aligned}\right.
\end{equation}
where $f(a)=-\frac{a}{1+a}$.

First, we establish the global well-posedness of the strong solution to the Cauchy problem \eqref{A4} within the critical regularity framework, as detailed below.

\begin{thm}[Global well-posedness]\label{T1.2}
   Let $d\geq 2$. There exists a positive constant $\delta_0>0$ such that if the initial data $(a_0,u_0)$ satisfy 
$a_0\in \dot B_{2,1}^{\frac{d}{2}}$, $u_0\in \dot B_{2,1}^{\frac{d}{2}-1}$, and 
\begin{align}\label{TB.1}
\|a_0\|_{\dot B_{2,1}^\frac{d}{2}}+   \|u_0\|_{\dot B_{2,1}^{\frac{d}{2}-1}} \leq \delta_0,
\end{align} 
then the   Cauchy problem \eqref{A4} admits a unique global strong solution 
$(a,u)$  satisfying
\begin{align}\label{TB.2}
a\in \mathcal{C}(\mathbb R^+; \dot B_{2,1}^{\frac{d}{2}}) \,\,\, \text{and}\,\,\,
u \in    \mathcal{C}(\mathbb R^+;  \dot B_{2,1}^{\frac{d}{2}-1})\cap L^1(\mathbb R^+;\dot B_{2,1}^{\frac{d}{2}+1}).          
\end{align}
Moreover, there exists a positive constant $C_1$ in dependent of the time $t$ 
such that, for any $t>0$,
\begin{align}\label{TB.3}
\|a\|_{\widetilde L_t^\infty( \dot B_{2,1}^{\frac{d}{2}})}+\|u\|_{\widetilde L^\infty_t( \dot B_{2,1}^{\frac{d}{2}-1})}+\|u\|_{L^1_t( \dot B_{2,1}^{\frac{d}{2}+1})}\leq C_1 \Big(\|a_0\|_{\dot B_{2,1}^\frac{d}{2}}+\|u_0\|_{\dot B_{2,1}^{\frac{d}{2}-1}}\Big).
\end{align}
\end{thm}

\begin{rem}
Compared with \cite{GTZ-JDE-2025}, where the condition  
\begin{align*}  
\|a_0\|_{H^3} + \|u_0\|_{H^4}+\|u_0\|_{L^1} \leq \delta_0 
\end{align*}  
is required, our assumption \eqref{TB.1} demands lower regularity. Moreover, we establish global well-posedness in the weaker critical Besov space $\dot B_{2,1}^{\frac{d}{2}} \times \dot B_{2,1}^{\frac{d}{2}-1}$ without requiring the initial data $u_0$ to be small in $L^1$, further highlighting the generality of our result.
\end{rem}

Next, we show the uniform stability of the pressureless Navier–Stokes system \eqref{A4}$_1$--\eqref{A4}$_2$.

\begin{thm}[Uniform stability]\label{T1.3}
 Let $d\geq 2$ and $\delta_0$ be the constant given in Theorem \ref{T1.2}.   
There
exists a small constant  $\delta_1\in(0,\delta_0)$, such that if
\begin{align*}
\max\Big\{\|a_0\|_{\dot B_{2,1}^\frac{d}{2} }+\|u_0\|_{\dot B_{2,1}^{\frac{d}{2}-1} }, \|\bar a_0\|_{\dot B_{2,1}^\frac{d}{2} }+\|\bar u_0\|_{\dot B_{2,1}^{\frac{d}{2}-1} } \Big\}   \leq \delta_1, 
\end{align*}
then the solutions $ (a,u)$ and $(\bar  a,\bar  u)$ obtained in Theorem \ref{T1.2}, corresponding to 
the initial data $ (a_0,u_0)$ and $(\bar  a_0,\bar  u_0)$,  satisfy
\begin{align}\label{error}
&\|\widetilde a(t)\|_{\widetilde L_t^\infty(\dot B_{2,1}^\frac{d}{2} )}+ \|\widetilde u\|_{\widetilde L^\infty_t( \dot B_{2,1}^{\frac{d}{2}-1})}+\|\widetilde u\|_{L^1_t( \dot B_{2,1}^{\frac{d}{2}+1})}   \leq C_2\Big(\|\widetilde a_0\|_{\dot B_{2,1}^\frac{d}{2}}+\|\widetilde u_0\|_{\dot B_{2,1}^{\frac{d}{2}-1}}\Big),   
\end{align}
for any $t>0  $, where
\begin{align*}
 (\widetilde a ,\widetilde u ):=  ((a-\bar a),  (u-\bar u)),
\end{align*}
and $C_2$ is a positive constant independent of the time $t$.
\end{thm}

We now present the optimal time decay rate for the pressureless Navier-Stokes system \eqref{A4}$_1$--\eqref{A4}$_2$ in Theorems \ref{T1.4}--\ref{T1.5}.

\begin{thm}[Upper-bound: the bounded condition of low frequency case] \label{T1.4}
Let $d\geq 2$. Under the assumptions of Theorem \ref{T1.2}, if the initial datum $u_0$ additionally satisfies  
\begin{align}\label{TC.1}
u_0^\ell\in   \dot B_{2,\infty}^{\sigma_0} \quad \text{with}\quad \sigma_0\in \Big[-\frac{d}{2},\frac{d}{2}-1\Big),
\end{align}
then for all $t \geq 1$, there exists a universal constant $C_3 > 0$ such that  
\begin{align}\label{TC.2}
\|u(t)\|_{\dot B_{2,1}^\sigma}\leq&\, C_3 \delta_{*} (1+t)^{-\frac{1}{2}(\sigma-\sigma_0)}, \\  \label{TC.4}
\|u(t)\|^h_{\dot B_{2,1}^{\frac{d}{2}+1}}\leq&\, C_3 \delta_{*} (1+t)^{-\frac{1}{2}(\frac{d}{2}-1-\sigma_0)}, 
\end{align}
for all $\sigma\in(\sigma_0,\frac{d}{2}-1]$. Here $\delta_{*}$ is defined as
\begin{align}\label{TC.3}
\delta_{*}:=\|u_0^\ell\|_{\dot B_{2,\infty}^{\sigma_0}} +\|u_0^h\|_{\dot B_{2,1}^{\frac{d}{2}-1}} +\|a_0\|_{\dot B_{2,1}^\frac{d}{2}}.  
\end{align}
\end{thm}

\begin{rem}
As mentioned in  \cite[Abstract]{GTZ-JDE-2025}, the classical perturbation theory from \cite{MN-JMKU-1980} does not apply to this model due to the absence of pressure. Moreover, the nonlinear terms involving  the factor $\frac{1}{1+a}$, such as $\frac{a \nabla{\rm div}\,u}{1+a}$ and $\frac{a\Delta u}{1+a}$, 
  lead to derivative losses and thereby present significant challenges in establishing optimal time decay rates for $u$. This suggests that the coupling between the density $a$ and $u$ currently induces detrimental effects rather than beneficial ones. Indeed, the density   lacks a dissipative structure  which inhibits its decay. Consequently, it remains bounded as stated in Theorem \ref{T1.2}, and no decay rate is derived for it in this setting.  This behavior stands in sharp contrast to that observed in the classical isentropic compressible Navier-Stokes equations, where the density  typically exhibits decay estimates similar to those of $u$. 
\end{rem}

\begin{rem}
In particular, when $d=3$ and $\sigma_0 = -\frac{3}{2}$, by utilizing the embeddings $L^1(\mathbb{R}^3) \hookrightarrow \dot{B}_{2,1}^{-\frac{3}{2}}(\mathbb{R}^3)$ and $\dot{B}_{2,1}^{0}(\mathbb{R}^3) \hookrightarrow L^2(\mathbb{R}^3)$, we derive the following optimal time-decay estimate for $u$, in the sense of Matsumura and Nishida \cite{MN-PJASAMS-1979}:  
\begin{align*}
\|u(t)\|_{L^2} \lesssim (1+t)^{-\frac{3}{4}}.
\end{align*}
This result further illustrates that the low-frequency assumption in $\dot B_{2,\infty}^{\sigma_0}$ is strictly weaker than the classical $L^1$ smallness condition employed in \cite{GTZ-JDE-2025}, thereby highlighting its broader applicability.
\end{rem}

If the initial data $u_0^\ell$ is sufficiently small in $\dot B_{2,\infty}^{\sigma_0}$ and $a_0^\ell \in \dot B_{2,\infty}^{\sigma_0+1}$ for $\sigma_0 \in [-\frac{d}{2}, \frac{d}{2}-1)$,
then improved decay estimates can be achieved compared to those in \eqref{TC.2}--\eqref{TC.4}. Moreover, the optimal time decay rates for $\|u(t)\|_{\dot B_{2,1}^{\sigma}}$ can be extended from $\sigma\in(\sigma_0,\frac{d}{2}-1]$ to the broader range $\sigma\in(\sigma_0,\frac{d}{2}+1]$.

\begin{thm}[Upper-bound: the smallness condition of low frequency case] \label{T1.1}
Let $d\geq 2$. Under the assumptions of Theorem \ref{T1.2}, if  the initial data $(a_0,u_0)$ further satisfies  
\begin{align}\label{TE.1}
a^\ell_0\in \dot B_{2,\infty}^{\sigma_0+1},\quad\|u_0^\ell\|_{\dot B_{2,\infty}^{\sigma_0}}\leq\eps_1 \quad \text{with}\quad \sigma_0\in \Big[-\frac{d}{2},\frac{d}{2}-1\Big),
\end{align}
where $\eps_1$ is a sufficiently small positive constant,
then for all $t \geq 1$, there exists a universal constant $C_4> 0$ such that  
\begin{align}\label{TE.2}
\|u(t)\|_{\dot B_{2,1}^\sigma}\leq&\, C_4 \delta_{*} (1+t)^{-\frac{1}{2}(\sigma-\sigma_0)},    \\ \label{TE.4}
\|u(t)\|_{\dot B_{2,1}^{\frac{d}{2}+1}}^h\leq&\, C_4 \delta_{*} (1+t)^{-\frac{1}{2}(\frac{d}{2}+1-\sigma_0)}, 
\end{align}
for all $\sigma\in(\sigma_0,\frac{d}{2}+1]$. Here $\delta_{*}$ is defined as same as in \eqref{TC.3}.
\end{thm}

\begin{rem}
The additional condition $a_0^\ell\in \dot B_{2,\infty}^{\sigma_0+1}$ in \eqref{TE.1} is necessary indeed. Since $a$ lacks a dissipative structure, the coupling terms between $a$ and $u$, such as $f(a)\Delta u$ and $f(a)\nabla{\rm div}\, u$, introduce a loss of derivatives on the density $a$ \emph{(}see Remark \ref{NJKR4.1}\emph{)}. To compensate for this loss, it is essential to impose extra low-frequency regularity on it. This phenomenon differs from that in classical isentropic compressible Navier–Stokes equations, such as  \cite{X-X-2021-JDE,Danchin-IM-00}, where the density possesses a dissipative structure. 
\end{rem}

\begin{rem}
As pointed out in \cite[Section 5]{LSZ-arXiv-2025} on the coupled pressureless Euler/Navier-Stokes system, the authors require the condition that $\rho_0^\ell \in \dot B_{2,\infty}^{\sigma_0}$ for the mass equation $\partial_t \rho + \mathrm{div}\,(\rho u) = 0$. However, our equation for $a$ read as $\partial_t a + \mathrm{div}\,u = -\mathrm{div}\,(a u)$, which contains an additional linear term $\mathrm{div}\,u$. This structural difference necessitates the assumption that $a_0^\ell \in \dot B_{2,1}^{\sigma_0+1}$, in stark contrast to the setting in \cite[Section 5]{LSZ-arXiv-2025}.
\end{rem}

\begin{rem}
When $d = 3$,  Wang et al.  \cite{WWX-2025-arXiv} considers the case in which the initial density $\rho_0$ is bounded both from above and below, their assumptions on $\rho_0$ are weaker than ours. Nevertheless, they only obtain a decay rate of order $1.5$ for $u$, whereas we achieve a higher decay rate of order $2.5$ under our framework.
\end{rem}

\begin{rem}
Unlike the results in \cite{LS-SIMA-2023}, where $u^h$ exhibits faster decay in $\dot B_{2,1}^{\frac{d}{2}+1}$, due to the absence of derivative loss in the coupling term, the decay rate in \eqref{TE.4} does not accelerate in our case. This is because the density $a$ lacks a dissipative structure and therefore cannot contribute additional decay in our system.
\end{rem}

\begin{rem}
 For $p \geq 2$ and $t \geq 1$, denoting $\Lambda := (-\Delta)^{-\frac{1}{2}}$ and combining \eqref{TE.2} with the embedding $\dot B_{2,1}^{\frac{d}{2}-\frac{d}{p}}(\mathbb R^d) \hookrightarrow L^p(\mathbb R^d)$,
we obtain the following $L^p$-type time decay estimate for $u$:
\begin{align*}
\|\Lambda^\sigma u\|_{L^p}\lesssim (1+t)^{-\frac{1}{2}(\sigma+\frac{d}{2}-\frac{d}{p}-\sigma_0)},\quad \text{with}\quad   \sigma+\frac{d}{2}-\frac{d}{p}\in\Big(\sigma_0,\frac{d}{2}+1\Big].  
\end{align*}

\end{rem}

To demonstrate the optimality of the time decay rates of the strong solution $u$ in Theorems \ref{T1.4}--\ref{T1.1}, it is necessary to establish a lower bound of the time decay rate for $u$   in Theorem \ref{T1.5} below. Prior to this, we introduce a subset $\dot{\mathcal{B}}_{2,\infty}^{\sigma_1}$ of the Besov spaces $\dot B_{2,\infty}^{\sigma_1}$ with $\sigma_1 \in \mathbb{R}$ (see \cite[Section 3]{Bl-2016-SIMA}):
\begin{equation}\label{G1.12}
\dot{\mathcal{B}}_{2,\infty}^{\sigma_1}:=
\left\{ f\in\dot B_{2,\infty}^{\sigma_1}\,\bigg|
 \begin{array}{l}
 \exists \,  c_0,   M_0>0, \,\exists \, \{k_j\}_{j\in\mathbb N}\subset \mathbb Z,\,  \textrm{such that}\, \, k_j\rightarrow-\infty,  \\
  |k_j-k_{j-1}|\leq M_0, \,\,\textrm{and}\,\, 2^{\sigma_1k_j} \|\dot\Delta_{k_j}f\|_{L^2}\geq c_0
 \end{array}
   \right\}.
\end{equation}

\begin{thm}[Lower-bound] \label{T1.5}
   Let $d\geq 2$. Under the assumptions of Theorem \ref{T1.2}, if the initial data $(a_0,u_0)$ further satisfies
\begin{align*}
{a^\ell_0\in \dot {{B}}_{2,\infty}^{\sigma_0+1},\quad u_0^\ell\in \dot{\mathcal{B}}_{2,\infty}^{\sigma_0},\quad\|u_0^\ell\|_{\dot {{B}}_{2,\infty}^{\sigma_0}}\leq\eps_2 \quad \text{with}\quad \sigma_0\in[-\frac{d}{2},\frac{d}{2}-1)},
\end{align*}
where $\eps_2$  is a sufficiently small positive constant,
then for all $t \geq 1$, there exists two universal constants $c_5>0$ and  $C_5 > 0$ such that  
\begin{align}\label{TD.1}
c_5(1+t)^{-\frac{1}{2}(\sigma-\sigma_0)}\leq \|u(t)\|_{\dot B_{2,1}^\sigma}\leq C_5 (1+t)^{-\frac{1}{2}(\sigma-\sigma_0)},    
\end{align}
for all $\sigma\in(\sigma_0,\frac{d}{2}+1]$, where  $\dot{\mathcal{B}}_{2,\infty}^{\sigma_0}$ is defined in \eqref{G1.12}.
\end{thm}

\begin{rem}
 In \cite{GTZ-JDE-2025}, an additional assumption ${\rm div}\,u_0=0$ is required. In contrast, by exploiting the fact that the velocity equation reduces to a heat equation and by utilizing the orthogonality between the operators $\mathbb{P}$ and $\mathbb{Q}$, which decompose vector fields into divergence-free and potential components, respectively, we show that the condition ${\rm div}\,u_0=0$ is not necessary. 
\end{rem}

\begin{rem}
When $d = 3$, under even weaker regularity assumptions, we have achieved the optimal decay rate for $u$ with $2.5$-order regularity.  This result improves the optimal decay rate corresponding to $2$-order regularity previously established in \cite{GTZ-JDE-2025} in the framework of Sobolev space.   
\end{rem}

\subsection{Strategies in our proofs}
In the proof of Theorem \ref{T1.2}, a fundamental challenge in studying the global well-posedness for the compressible pressureless Navier-Stokes equations \eqref{A4} arises from the $L^1$ time integrability of $\nabla u$, namely:
\begin{align}\label{control-u}
\int_0^\infty \|\nabla u(t)\|_{L^\infty} \, {\rm d}t < \infty,
\end{align}
which is crucial for controlling both the linear term ${\rm div}\,u$ and the nonlinear term $a{\rm div}\,u$, and thereby closing the estimate for $a$. We note that in \cite{GTZ-JDE-2025}, to overcome the difficulty related to \eqref{control-u}, an additional smallness assumption on the initial datum $u_0$ in $L^1$ was imposed, which yields a decay rate of $\|\nabla u(t)\|_{H^2}$ as $(1+t)^{-\frac{5}{4}}$. Then, by using the embedding $H^2(\mathbb{R}^3) \hookrightarrow L^{\infty}(\mathbb{R}^3)$, it follows that
\begin{align*}
\int_0^\infty \|\nabla u(t)\|_{L^\infty} \, {\rm d}t \lesssim \int_0^\infty \|\nabla u(t)\|_{H^2} \, {\rm d}t < \infty.
\end{align*}
In \cite{WWX-2025-arXiv}, due to the absence of an estimate for $\|u\|_{L^1(\dot B_{2,1}^{\frac{d}{2}+1})}$, the challenge posed by \eqref{control-u} was addressed through estimates in Lorentz spaces.  Unlike \cite{GTZ-JDE-2025,WWX-2025-arXiv}, here we take advantage of the embedding $\dot B_{2,1}^{\frac{d}{2}}(\mathbb{R}^d)\hookrightarrow L^\infty(\mathbb{R}^d)$ to establish the $L^1$ time integrability of $\nabla u$  in the Besov space framework, specifically in $\dot B_{2,1}^{\frac{d}{2}}(\mathbb{R}^d)$.

When establishing the uniform stability of the system \eqref{A4}, we define the following functional (see \eqref{G4.2}):
\begin{align*} 
\widetilde {\mathcal{X}}(t):=  \|\widetilde  a\|_{\widetilde L_t^\infty( \dot B_{2,1}^{\frac{d}{2}})}+\|\widetilde  u\|_{\widetilde L^\infty_t( \dot B_{2,1}^{\frac{d}{2}-1})}+\|\widetilde u\|_{L^1_t( \dot B_{2,1}^{\frac{d}{2}+1})}.   
\end{align*}
Our goal is to establish the estimate:
 \begin{align*}
  \widetilde {\mathcal{X}}(t)\lesssim   \|\widetilde a_0\|_{\dot B_{2,1}^\frac{d}{2}}+\|\widetilde u_0\|_{\dot B_{2,1}^{\frac{d}{2}-1}}+(\delta_0+\delta_1)    \widetilde {\mathcal{X}}(t).
 \end{align*}
Although the estimate for $\widetilde u$ follows directly from Lemma \ref{LA.6}, the corresponding estimate for $\widetilde a$ cannot be derived in the same manner from Lemma \ref{LA.7}, due to the presence of the term $\|\widetilde u \cdot \nabla \bar a\|_{L_t^1(\dot B_{2,1}^{\frac{d}{2}})}$, which requires additional justification. As noted in Remark \ref{REM3.1}, this term is not easily controlled. However, this difficulty can be circumented by utilizing a more refined commutator estimate (see the proof in Lemma \ref{L4.2}).

To obtain the optimal time decay rate of $u$ for the problem \eqref{A4}, we establish the upper and lower bounds of the estimates for $u$. For the upper bound, we first consider the case under the assumption that $u^\ell \in \dot B_{2,\infty}^{\sigma_0}$ with $\sigma_0 \in [-\frac{d}{2}, \frac{d}{2}-1)$. We first prove the propagation of $\dot B_{2,\infty}^{\sigma_0}$ in the low-frequency regime in Lemma \ref{L5.1}. Then, using the method in \cite{LS-SIMA-2023}, and by establishing time-weighted estimates at both high and low frequencies in Lemmas \ref{L5.1}--\ref{L5.2}, we derive an estimate for $\mathcal X_{M}(t)$ (see \eqref{G5.23}):
\begin{align*}
\mathcal{X}_{M}(t) := \|\tau^M u\|_{\widetilde L_t^\infty(\dot B_{2,1}^{{\frac{d}{2}}-1})} + \|\tau^M u\|_{L_t^1(\dot B_{2,1}^{\frac{d}{2}+1})} \lesssim \delta_{*} t^{M - \frac{1}{2}(\frac{d}{2}-1-\sigma_0)}.
\end{align*}
Thus, \eqref{TC.2} follows. Furthermore, \eqref{TC.4} is obtained via maximal regularity estimates provided in Lemma \ref{LA.6}. Hence, the proof of Theorem \ref{T1.4} is completed. In fact, the upper bound of the optimal time decay rates for $\|u(t)\|_{\dot B_{2,1}^{\sigma}}$, initially established for $\sigma\in(\sigma_0,\frac{d}{2}-1]$ in Theorem \ref{T1.4}, can be extended to the broader range $\sigma\in(\sigma_0,\frac{d}{2}+1]$. To achieve this extension, we further assume that $u_0 \in \dot B_{2,\infty}^{\sigma_0}$ is sufficiently small and $a_0^\ell \in \dot B_{2,\infty}^{\sigma_0+1}$, where $\sigma_0 \in [-\frac{d}{2}, \frac{d}{2}-1)$. The assumption on $a_0^\ell$ is technical in nature and is introduced primarily to handle the derivative loss arising in the coupled terms $f(a)\Delta u$ and $f(a)\nabla {\rm div}\,u$ (see Remark \ref{NJKR4.1}). This difficulty is fundamentally different from that encountered in the isentropic compressible Navier-Stokes equations; in contrast, as demonstrated in \cite{Danchin-IM-00}, the density variable $a$ possesses both a dissipative structure and a well-defined time-decay rate. Our approach in Theorem \ref{T1.1} is largely inspired by \cite[Section 5.2]{Danchin-2018}.

Finally, to derive the lower bound estimate for $u$, we decompose the solution $u$ into a linear component $u_L$ and a nonlinear component $\omega$. For the linearized system \eqref{G5.26}, by adapting the approach employed in \cite{Bl-2016-SIMA} and \cite[Section 3]{BSXZ-Adv-2024}, we establish the linear analysis of $u_L$ in Lemma \ref{L5.4}. The nonlinear analysis for $\omega$, provided in Lemma \ref{L4.7}, follows an argument analogous to that in Lemma \ref{NJKL4.4}. By applying Duhamel’s principle and combining the estimates from both the linear and nonlinear components, we complete the proof of Theorem \ref{T1.5}.

\subsection{Structure of our paper}
The remainder of this paper is structured as follows. In Section 2, we derive the a priori estimates for the solution   $(a,u)$ and establish the global existence of solutions to the pressureless Navier–Stokes system \eqref{A4}. Section 3 is devoted to proving the uniform stability of the strong solution $(a,u)$. In Section 4, we demonstrate the optimal decay rate of  $u$ in $L^2$ norm. For the upper bound estimate, we consider two distinct cases: one under a low-frequency bounded initial condition and the other under a smallness assumption; we show that the latter leads to a higher-order decay rate. For the lower bound, we establish the convergence rate of  $u$ through a combination of linear and nonlinear analysis. Finally, Appendix A 
collects several fundamental properties of Besov spaces and product estimates, which have been frequently used in previous sections.

\section{Global well-posedness  of the pressureless Navier–Stokes system}
In this section, we prove Theorem \ref{T1.2}, which concerns the global existence and uniqueness of the solution $(a,u)$ to the Cauchy problem \eqref{A4}. We first establish the following  uniform-in-time {\it a priori estimates}.
\begin{prop}\label{P3.1}
Assume that $(a,u)$ is a strong solution to
 the Cauchy problem \eqref{A4} defined on $[0,T)\times \mathbb R^d$ with a given time $T>0$.  For any $t\in (0,T)$, it holds that, for a given generic constant $0<\delta<1$ to be chosen later, if the solution 
 $(a,u)$ satisfy 
\begin{align}\label{G3.1}
 \|a\|_{\widetilde L_t^\infty(\dot B_{2,1}^{\frac{d}{2}})}+\|u\|_{\widetilde L_t^\infty(\dot B_{2,1}^{\frac{d}{2}-1})}  \leq \delta,
\end{align}
then the following estimates holds 
\begin{align}\label{G3.2}
 \|a\|_{\widetilde L_t^\infty(\dot B_{2,1}^{\frac{d}{2}})}+\|u\|_{\widetilde L_t^\infty(\dot B_{2,1}^{\frac{d}{2}-1})}+\|u\|_{  L_t^1(\dot B_{2,1}^{\frac{d}{2}+1})}\leq C_0\Big(\|a_0\|_{ \dot B_{2,1}^{\frac{d}{2}}}+\|u_0\|_{\dot B_{2,1}^{\frac{d}{2}-1}}\Big),   
\end{align}
where $C_0>0$ is a constant independent of the time $T$.
\end{prop}
The proof of Proposition \ref{P3.1} is based on Lemmas \ref{L3.2}--\ref{L3.4} below.

\subsection{Estimate of $u$}
We begin by estimating $u$ in both the low-frequency and high-frequency components.
\begin{lem}\label{L3.2}
Let $(a,u)$ be a strong solution to the Cauchy problem \eqref{A4} on $[0,T)\times \mathbb R^d$. Then, under the condition \eqref{G3.1}, it holds that
\begin{align}\label{G3.4}
\|u\|_{\widetilde L_t^\infty(\dot B_{2,1}^{\frac{d}{2}-1})}^\ell+ \|u\|_{  L_t^1(\dot B_{2,1}^{\frac{d}{2}+1})}^\ell\lesssim  \|u_0\|_{\dot B_{2,1}^{\frac{d}{2}-1}}+\delta    \|u\|_{  L_t^1(\dot B_{2,1}^{\frac{d}{2}+1})}.
\end{align}
\end{lem}
\begin{proof}
Taking the $L^2$ inner product of \eqref{A4}$_2$ with $\dot\Delta_k u$, we arrive at 
\begin{align}\label{G3.5}
&\frac{1}{2}\frac{{\rm d}}{{\rm d}t}\|\dot\Delta_k u\|_{L^2}^2+\mu \|\nabla \dot\Delta_k u\|_{L^2}^2+(\mu+\nu)\|{\rm div}\,\dot\Delta_k u\|_{L^2}^2  \nonumber\\
&\quad \lesssim  \|\dot\Delta_k u\|_{L^2}\big(\|\dot\Delta_k(u\cdot\nabla u)\|_{L^2}+\|\dot\Delta_k(f(a)\Delta u)\|_{L^2}+\|\dot\Delta_k(f(a)\nabla{\rm div}\, u)\|_{L^2}        \big), 
\end{align}
for $k\leq 0$.
Dividing \eqref{G3.5} by $\big({\|\dot\Delta_k u\|^2_{L^2}} + \eps_{*}^2\big)^{\frac{1}{2}}$ with $\eps_{*} > 0$, we integrate the resulting inequality over $[0,t]$ and then pass to the limit as $\eps_{*} \to 0$ to obtain
\begin{align}\label{G3.6}
&\|\dot\Delta_k u\|_{L^2} +2^{2k}\int_0^t \|\dot\Delta_k u\|_{L^2}{\rm d}\tau\nonumber\\
&\quad \lesssim  \|\dot\Delta_k u_0\|_{L^2} +\int_0^t  \big(\|\dot\Delta_k(u\cdot\nabla u)\|_{L^2}+\|\dot\Delta_k(f(a)\Delta u)\|_{L^2}+\|\dot\Delta_k(f(a)\nabla{\rm div}\, u)\|_{L^2}        \big){\rm d}\tau.
\end{align}
Multiplying \eqref{G3.6} by $2^{k(\frac{d}{2}-1)}$, taking the supremum over $[0,t]$, and summing over all $k \leq 0$,  one gets 
\begin{align}\label{G3.7}
&\|u\|_{\widetilde L_t^\infty(\dot B_{2,1}^{\frac{d}{2}-1})}^\ell+\|u\|_{L_t^1(\dot B_{2,1}^{\frac{d}{2}+1})}^\ell\nonumber\\
&\quad \lesssim  \|u_0\|_{ \dot B_{2,1}^{\frac{d}{2}-1}}^\ell +\|u\cdot\nabla u\|_{L_t^1(\dot B_{2,1}^{\frac{d}{2}-1})}^\ell+ \|f(a)\Delta u\|_{L_t^1(\dot B_{2,1}^{\frac{d}{2}-1})}^\ell+  \|f(a)\nabla{\rm div}\, u\|_{L_t^1(\dot B_{2,1}^{\frac{d}{2}-1})}^\ell.
\end{align}
From the estimate \eqref{A.3} in Lemma \ref{LA.3}, it follows that
\begin{align}\label{G3.8}
\|u\cdot\nabla u\|_{L_t^1(\dot B_{2,1}^{\frac{d}{2}-1})} \lesssim   \|u\|_{L_t^\infty(\dot B_{2,1}^{\frac{d}{2}-1})}\|u\|_{L_t^1(\dot B_{2,1}^{\frac{d}{2}+1})}\lesssim \delta \|u\|_{L_t^1(\dot B_{2,1}^{\frac{d}{2}+1})}.
\end{align}
By leveraging Lemma \ref{LA.3} again and the continuity of composition functions  in Lemma \ref{LA.4},
 we derive  that 
\begin{align}\label{G3.9}
 \|f(a)\Delta u\|_{L_t^1(\dot B_{2,1}^{\frac{d}{2}-1})} +  \|f(a)\nabla{\rm div}\, u\|_{L_t^1(\dot B_{2,1}^{\frac{d}{2}-1})}   \lesssim&\, \|f(a)\|_{L_t^\infty(\dot B_{2,1}^\frac{d}{2})}\|u\|_{L_t^1(\dot B_{2,1}^{\frac{d}{2}+1})} \nonumber\\
 \lesssim&\, \|a\|_{\widetilde L_t^\infty(\dot B_{2,1}^\frac{d}{2})}\|u\|_{L_t^1(\dot B_{2,1}^{\frac{d}{2}+1})}\nonumber \\
 \lesssim&\, \delta\|u\|_{L_t^1(\dot B_{2,1}^{\frac{d}{2}+1})}.
\end{align}
Putting the estimates \eqref{G3.8} and \eqref{G3.9} into \eqref{G3.7} gives rise to \eqref{G3.4}.
\end{proof}

\begin{lem}\label{L3.3}
Let $(a,u)$ be a strong solution to the Cauchy problem \eqref{A4} on $[0,T)\times \mathbb R^d$. Then, under the condition \eqref{G3.1}, it holds that
\begin{align}\label{G3.10}
\|u\|_{\widetilde L_t^\infty(\dot B_{2,1}^{\frac{d}{2}-1})}^h+ \|u\|_{  L_t^1(\dot B_{2,1}^{\frac{d}{2}+1})}^h\lesssim  \|u_0\|_{\dot B_{2,1}^{\frac{d}{2}-1}}+\delta    \|u\|_{  L_t^1(\dot B_{2,1}^{\frac{d}{2}+1})}.
\end{align}
\end{lem}
\begin{proof}
Noting that $\dot\Delta_k(u\cdot\nabla u) = u\cdot\nabla\dot\Delta_k u - [u\cdot\nabla,\dot\Delta_k]u$, it follows from \eqref{A4}$_2$ that
\begin{align}\label{G3.11}
&\frac{1}{2}\frac{{\rm d}}{{\rm d}t}\|\dot\Delta_k u\|_{L^2}^2+\mu \|\nabla \dot\Delta_k u\|_{L^2}^2+(\mu+\nu)\|{\rm div}\,\dot\Delta_k u\|_{L^2}^2  \nonumber\\
&\quad \lesssim  \|\dot\Delta_k u\|_{L^2}\big(\|[u\cdot\nabla,\dot\Delta_k] u\|_{L^2}+\|\dot\Delta_k(f(a)\Delta u)\|_{L^2} +\|\dot\Delta_k(f(a)\nabla{\rm div}\, u)\|_{L^2}        \big)\nonumber\\
&\qquad+\|{\rm div}\, u\|_{L^\infty}\|\dot\Delta_k u\|_{L^2}^2,
\end{align}
for $k\geq-1$,  which implies that 
\begin{align}\label{G3.12}
&\|\dot\Delta_k u\|_{L^2} +2^{2k}\int_0^t \|\dot\Delta_k u\|_{L^2}{\rm d}\tau\nonumber\\
&\quad \lesssim  \|\dot\Delta_k u_0\|_{L^2} +\int_0^t  \big(\|{\rm div}\, u\|_{L^\infty}\|\dot\Delta_k u\|_{L^2}+\|[u\cdot\nabla,\dot\Delta_k] u\|_{L^2} +\|\dot\Delta_k(f(a)\Delta u)\|_{L^2}    \big){\rm d}\tau\nonumber\\
&\qquad+\int_0^t\|\dot\Delta_k(f(a)\nabla{\rm div}\, u)\|_{L^2}  {\rm d}\tau.
\end{align}
A direct computation for high frequencies ($k\geq -1$) yields  
\begin{align}\label{G3.13}
&\|u\|_{\widetilde L_t^\infty(\dot B_{2,1}^{\frac{d}{2}-1})}^h+\|u\|_{L_t^1(\dot B_{2,1}^{\frac{d}{2}+1})}^h\nonumber\\
&\quad \lesssim  \|u_0\|_{ \dot B_{2,1}^{\frac{d}{2}-1}}^h +\|{\rm div}\,u\|_{L_t^1(L^\infty)}\|u \|_{\widetilde L_t^\infty(\dot B_{2,1}^{\frac{d}{2}-1})}^h+\sum_{k\geq -1}2^{k(\frac{d}{2}-1)}\|[u\cdot\nabla,\dot\Delta_k] u\|_{L_t^1(L^2)} \nonumber\\
&\qquad+ \|f(a)\Delta u\|_{L_t^1(\dot B_{2,1}^{\frac{d}{2}-1})}^h+  \|f(a)\nabla{\rm div}\, u\|_{L_t^1(\dot B_{2,1}^{\frac{d}{2}-1})}^h.  
\end{align}
Applying the embedding $\dot B_{2,1}^\frac{d}{2}(\mathbb{R}^d) \hookrightarrow L^\infty(\mathbb{R}^d)$ in Lemma \ref{LA.2} yields  
\begin{align}\label{G3.14}
\|{\rm div}\,u\|_{L_t^1(L^\infty)}\|u \|_{\widetilde L_t^\infty(\dot B_{2,1}^{\frac{d}{2}-1})}^h\lesssim \|u\|_{L_t^1(\dot B_{2,1}^{\frac{d}{2}+1})}\|u\|_{\widetilde L_t^\infty(\dot B_{2,1}^{\frac{d}{2}-1})} \lesssim\delta \|u\|_{L_t^1(\dot B_{2,1}^{\frac{d}{2}+1})}.  
\end{align}
According to the commutator estimate \eqref{A.6} established in Lemma \ref{LA.5}, we have
\begin{align}\label{G3.15}
 \sum_{k\geq -1}2^{k(\frac{d}{2}-1)}\|[u\cdot\nabla,\dot\Delta_k] u\|_{L_t^1(L^2)}\lesssim \|u\|_{L_t^1(\dot B_{2,1}^{\frac{d}{2}+1})} \|u\|_{\widetilde L_t^\infty(\dot B_{2,1}^{\frac{d}{2}-1})}\lesssim \delta \|u\|_{L_t^1(\dot B_{2,1}^{\frac{d}{2}+1})}.  
\end{align}
Inserting the estimates \eqref{G3.9} and \eqref{G3.14}--\eqref{G3.15} into \eqref{G3.13}, we consequently obtain the desired estimate \eqref{G3.10}.
\end{proof}
Combining Lemmas \ref{L3.2} and \ref{L3.3} up, we directly get
\begin{align}\label{G3.16}
\|u\|_{\widetilde L_t^\infty(\dot B_{2,1}^{\frac{d}{2}-1})} + \|u\|_{  L_t^1(\dot B_{2,1}^{\frac{d}{2}+1})} \lesssim  \|u_0\|_{\dot B_{2,1}^{\frac{d}{2}-1}}+\delta    \|u\|_{  L_t^1(\dot B_{2,1}^{\frac{d}{2}+1})}.
\end{align}

\subsection{Estimate of $a$}
Finally, by exploiting the transport structure of \eqref{A4}$_1$ and the Lipschitz bound on $u$, we establish the estimate of $a$.

\begin{lem}\label{L3.4}
Let $(a,u)$ be a strong solution to the Cauchy problem \eqref{A4} on $[0,T)\times \mathbb R^d$. Then, under the condition \eqref{G3.1}, it holds that
\begin{align}\label{G3.17}
\|a\|_{\widetilde L_t^\infty(\dot B_{2,1}^{\frac{d}{2}-1})} \lesssim   \|a_0\|_{\dot B_{2,1}^{\frac{d}{2}-1}}+  (1+\delta)\|u\|_{  L_t^1(\dot B_{2,1}^{\frac{d}{2}+1})}.
\end{align}
\end{lem}
\begin{proof}
Applying the operator $\dot\Delta_k$ to \eqref{A4}$_1$, one has
\begin{align}\label{G3.18}
\partial_t\dot\Delta_k a+{\rm div}\, \dot\Delta_k u+\dot\Delta_k(a{\rm div}\, u)+u\cdot\nabla\dot\Delta_k a-[u\cdot\nabla, \dot\Delta_k]a=0.
\end{align}
Thanks to the standard $L^2$ energy estimate on \eqref{G3.18}, we find that
\begin{align*}
 \frac{1}{2}\frac{{\rm d}}{{\rm d}t}\|\dot\Delta_k a\|_{L^2}^2=&\,\frac{1}{2}\int_{\mathbb R^d}{\rm div}\,u|\dot\Delta_k a|^2{\rm d}x+\int_{\mathbb R^d}[u\cdot\nabla, \dot\Delta_k]a\dot\Delta_k a{\rm d}x\nonumber\\
 &-\int_{\mathbb R^d} \dot\Delta_k(a{\rm div}\, u)\dot\Delta_k a{\rm d}x-\int_{\mathbb R^d} {\rm div}\, \dot\Delta_k u \dot\Delta_k a{\rm d}x \nonumber\\
  \lesssim&\, \big(\|[u\cdot\nabla, \dot\Delta_k]a\|_{L^2}+\|\dot\Delta_k(a{\rm div}\, u)\|_{L^2}+2^k\|\dot\Delta_ku\|_{L^2}\big)  \|\dot\Delta_k a\|_{L^2}\nonumber\\
&+\|{\rm div}\, u\|_{L^\infty}\|\dot\Delta_k a\|_{L^2}^2,
\end{align*}
which leads to
\begin{align}\label{G3.19}
 \|\dot\Delta_k a\|_{L^2}  
  \lesssim&\, \|\dot\Delta_k a_0\|_{L^2}+\int_0^t\big( \|[u\cdot\nabla, \dot\Delta_k]a\|_{L^2}+\|\dot\Delta_k(a{\rm div}\, u)\|_{L^2}+2^k\|\dot\Delta_ku\|_{L^2}\big){\rm d}\tau\nonumber\\
  &+\int_0^t\|{\rm div}\, u\|_{L^\infty}\|\dot\Delta_k a\|_{L^2}{\rm d}\tau.
\end{align}
Multiplying \eqref{G3.19} by $2^{\frac{kd}{2}}$, taking the supremum over $[0,t]$, and summing over all $k\in\mathbb{Z}$,  we arrive at 
\begin{align}\label{G3.20}
\|a\|_{\widetilde L_t^\infty(\dot B_{2,1}^\frac{d}{2})}\lesssim&\, \|a_0\|_{\dot B_{2,1}^\frac{d}{2}}+    \sum_{k\in \mathbb Z}2^{\frac{kd}{2}}\|[u\cdot\nabla,\dot\Delta_k] a\|_{L_t^1(L^2)}+\|a{\rm div}\, u\|_{L_t^1(\dot B_{2,1}^\frac{d}{2})}\nonumber\\
&+\|u\|_{L_t^1(\dot B_{2,1}^{\frac{d}{2}+1})}+\|{\rm div}\, u\|_{L_t^1( L^\infty)}\|a\|_{\widetilde L_t^\infty(\dot B_{2,1}^\frac{d}{2})}.
\end{align}
By virtue of the embedding $\dot B_{2,1}^\frac{d}{2}(\mathbb{R}^d)\hookrightarrow L^\infty(\mathbb{R}^d)$ and Lemmas \ref{LA.2}--\ref{LA.3} and \ref{LA.5}, we compute that
\begin{align*}
\sum_{k\in \mathbb Z}2^{\frac{kd}{2}}\|[u\cdot\nabla,\dot\Delta_k] a\|_{L_t^1(L^2)}\lesssim&\, \|u\|_{L_t^1(\dot B_{2,1}^{\frac{d}{2}+1})}\|a\|_{\widetilde L_t^\infty(\dot B_{2,1}^\frac{d}{2})}\lesssim \delta    \|u\|_{L_t^1(\dot B_{2,1}^{\frac{d}{2}+1})},\\
\|a{\rm div}\, u\|_{L_t^1(\dot B_{2,1}^\frac{d}{2})}\lesssim&\,  \|u\|_{L_t^1(\dot B_{2,1}^{\frac{d}{2}+1})}\|a\|_{\widetilde L_t^\infty(\dot B_{2,1}^\frac{d}{2})}\lesssim \delta    \|u\|_{L_t^1(\dot B_{2,1}^{\frac{d}{2}+1})},\\
\|{\rm div}\, u\|_{L_t^1( L^\infty)}\|a\|_{\widetilde L_t^\infty(\dot B_{2,1}^\frac{d}{2})}\lesssim&\, \|u\|_{L_t^1(\dot B_{2,1}^{\frac{d}{2}+1})}\|a\|_{\widetilde L_t^\infty(\dot B_{2,1}^\frac{d}{2})}\lesssim \delta    \|u\|_{L_t^1(\dot B_{2,1}^{\frac{d}{2}+1})}.
\end{align*}
By substituting the above estimates into \eqref{G3.20}, we derive \eqref{G3.17}.
\end{proof}

With the aid of Lemmas \ref{L3.2}--\ref{L3.4}, we now proceed to prove Proposition \ref{P3.1}.

\begin{proof}[Proof of Proposition \ref{P3.1}]
It follows from \eqref{G3.16} that  
\begin{align*}
\|u\|_{\widetilde L_t^\infty(\dot B_{2,1}^{\frac{d}{2}-1})} + \|u\|_{  L_t^1(\dot B_{2,1}^{\frac{d}{2}+1})} \leq C_* \|u_0\|_{\dot B_{2,1}^{\frac{d}{2}-1}}+C_*\delta    \|u\|_{  L_t^1(\dot B_{2,1}^{\frac{d}{2}+1})},    
\end{align*}
for some uniform constant $C_* > 0$. By selecting $\delta\leq\min\big\{1, \frac{1}{2C*}\big\}$, we obtain that
\begin{align}\label{G3.21}
 \|u\|_{\widetilde L_t^\infty(\dot B_{2,1}^{\frac{d}{2}-1})} + \|u\|_{  L_t^1(\dot B_{2,1}^{\frac{d}{2}+1})} \leq2 C_* \|u_0\|_{\dot B_{2,1}^{\frac{d}{2}-1}}.   
\end{align}
On the other hand, from \eqref{G3.17}, we have
\begin{align*}
\|a\|_{\widetilde L_t^\infty(\dot B_{2,1}^{\frac{d}{2}-1})} \leq C   \|a_0\|_{\dot B_{2,1}^{\frac{d}{2}-1}}+ C\|u\|_{  L_t^1(\dot B_{2,1}^{\frac{d}{2}+1})},    
\end{align*}
which together with \eqref{G3.21} gives
\begin{align}\label{G3.22}
\|a\|_{\widetilde L_t^\infty(\dot B_{2,1}^{\frac{d}{2}-1})} \leq C   \|a_0\|_{\dot B_{2,1}^{\frac{d}{2}-1}}+ 2CC_*\|u_0\|_{  \dot B_{2,1}^{\frac{d}{2}-1} }.  
\end{align}
By combining \eqref{G3.21} with \eqref{G3.22}, we directly obtain \eqref{G3.2}, thereby completing the proof of Proposition \ref{P3.1}.
\end{proof}

\subsection{global well-posedness of strong solution }
This subsection aims to prove Theorem \ref{T1.2}. To begin with, we state the existence and uniqueness of local-in-time solution for the Cauchy problem \eqref{A4}, which can be established by using a method similar to that in \cite[Section 3]{Danchin-2018}. For brevity, we omit the details here.

\begin{thm}[Local well-posedness]\label{T3.5}
  Let $d\geq 2$ and  assume that the initial data $(a_0,u_0)$ satisfies
\begin{align}\label{TA.1}
a_0\in \dot B_{2,1}^{\frac{d}{2}},\quad \inf_{x\in\mathbb R^d} (1+a_0)(x)>0, \quad u_0\in     \dot B_{2,1}^{\frac{d}{2}-1}.
\end{align}
Then, there exists a time $T>0$, such that the Cauchy problem \eqref{A4} admits a unique strong solution $(a,u)$ satisfying that for $t\in[0,T)$,
\begin{equation}\label{TA.2}\left\{
\begin{aligned}
&    a\in \mathcal{C}([0,T); \dot B_{2,1}^{\frac{d}{2}}),\quad \inf_{x\in\mathbb R^d} (1+a )(x,t)>0,      \\
&  u \in    \mathcal{C}([0,T); \dot B_{2,1}^{\frac{d}{2}-1})\cap L^1(0,T;\dot B_{2,1}^{\frac{d}{2}+1}).          
\end{aligned}\right.
\end{equation}
\end{thm}

\begin{proof}[Proof of Theorem \ref{T1.2}]
Theorem \ref{T3.5} guarantees the existence of a maximal existence time $T_0$, such that the Cauchy problem \eqref{A4} admits a unique strong solution $(a,u)$ satisfying \eqref{TA.2}. Define 
\begin{align}\label{defnx}
\mathcal{X}(t):=\|a\|_{\widetilde L_t^\infty(\dot B_{2,1}^\frac{d}{2})}+\|u\|_{\widetilde L_t^\infty(\dot B_{2,1}^{\frac{d}{2}-1})}+\|u\|_{  L_t^1(\dot B_{2,1}^{\frac{d}{2}+1})}.    
\end{align}
We then set 
\begin{align}\label{G3.24}
T^*:=\sup \{t\in[0,T_0)\, | \,  \mathcal{X}(t)\leq \delta   \},   
\end{align}
 and claim that $T^*=T_0$. Suppose, for contradiction, that $T^*<T_0$. Choosing
\begin{align*}
\delta_0:=\min\Big\{1,\frac{\delta}{2C_0}\Big\},    
\end{align*}
where $C_0$ is the constant given in Proposition \ref{P3.1} and $\delta_0$ is defined by \eqref{TB.1}. 
Then, by the a priori estimates \eqref{G3.2} in Proposition \ref{P3.1}, we have
\begin{align*}
\mathcal{X}(t)\leq C_0    \Big(\|a_0\|_{ \dot B_{2,1}^{\frac{d}{2}}}+\|u_0\|_{\dot B_{2,1}^{\frac{d}{2}-1}}\Big) \leq \frac{\delta}{2} ,
\end{align*}
for all $0<t<T^*$. Since $\mathcal{X}(t)$ is continuous in time, it follows that
\begin{align*}
\mathcal{X}(T^*) \leq \frac{\delta}{2} , 
\end{align*}
which contradicts the definition of $T^*$ in \eqref{G3.24}. Therefore, $T^*=T_0$  holds.

Finally, we claim that $T^* = T_0 = +\infty$. If $T_0 < +\infty$, then by combining Theorem \ref{T3.5} with the uniform estimate \eqref{G3.2} in Proposition \ref{P3.1}, the strong solution $(a,u)$ can be extended to an interval $[0, T_0 + \eta_0]$ for some constant $\eta_0 > 0$. This contradicts the definition of $T_0$ as the maximal existence time. Therefore, $(a,u)$ constitutes a global strong solution to the Cauchy problem \eqref{A4}, and the estimate \eqref{TB.3} holds. The proof of Theorem \ref{T1.2} is thus completed.
\end{proof}

\section{Uniform stability of the pressureless Navier–Stokes system }

This section studies the uniform stability of the unique strong solution established in Theorem \ref{T1.2}. Specifically, we consider two solutions, $(a, u)$ and $(\bar{a}, \bar{u})$, to the pressureless Navier--Stokes system \eqref{A4}$_1$--\eqref{A4}$_2$, corresponding to the initial data $(a_0, u_0)$ and $(\bar{a}_0, \bar{u}_0)$, respectively. Under the assumptions of Theorem \ref{T1.3}, both solutions satisfy the following inequalities for all $t > 0$:
\begin{align*}
\|a\|_{\widetilde L_t^\infty( \dot B_{2,1}^{\frac{d}{2}})}+\|u\|_{\widetilde L^\infty_t( \dot B_{2,1}^{\frac{d}{2}-1})}+\|u\|_{L^1_t( \dot B_{2,1}^{\frac{d}{2}+1})}\lesssim&\,   \|a_0\|_{\dot B_{2,1}^\frac{d}{2}}+\|u_0\|_{\dot B_{2,1}^{\frac{d}{2}-1}} ,      \\
\|\bar a\|_{\widetilde L_t^\infty( \dot B_{2,1}^{\frac{d}{2}})}+\|\bar u\|_{\widetilde L^\infty_t( \dot B_{2,1}^{\frac{d}{2}-1})}+\|\bar u\|_{L^1_t( \dot B_{2,1}^{\frac{d}{2}+1})}\lesssim&\,   \|\bar a_0\|_{\dot B_{2,1}^\frac{d}{2}}+\|\bar u_0\|_{\dot B_{2,1}^{\frac{d}{2}-1}} .
\end{align*}
To begin with, we examine the error equations between $(a, u)$ and $(\bar a,\bar u)$. 
Set  $(\widetilde a ,\widetilde u )=  ((a-\bar a),  (u-\bar u))$.  
From \eqref{A4}, it follows that  
\begin{equation}\label{G4.1}
\left\{
\begin{aligned}
&\partial_t\widetilde a+u\cdot\nabla \widetilde a=-\bar a{\rm div}\, \widetilde u-\widetilde a{\rm div}\, u-\widetilde u\cdot\nabla\bar a  -{\rm div}\, \widetilde u,\\
&\partial_t\widetilde{u}-\mu\Delta\widetilde u-(\mu+\nu)\nabla{\rm div}\, \widetilde u=-\widetilde u\cdot\nabla\bar u-u\cdot\nabla\widetilde u+\widetilde F_1+\widetilde F_2  ,
\end{aligned}\right.
\end{equation}
where $\widetilde F_1$ and $\widetilde F_2$ are given by
\begin{align*}
\widetilde F_1:=&\, \mu \big(f(a )-f(\bar a)\big)\Delta u +\mu f(\bar a)\Delta\widetilde u,  \\ 
\widetilde F_2:=&\, (\mu+\nu) \big(f(a )-f(\bar a )\big)\nabla{\rm div}\, u +(\mu+\nu) f(\bar a)\nabla{\rm div}\,\widetilde u,
\end{align*}
with the initial data
\begin{align*}
(\widetilde a,\widetilde u)|_{t=0}=(\widetilde a_0(x), \widetilde u_0(x)).    
\end{align*}

Next, we establish the $\dot B_{2,1}^\frac{d}{2}$-regularity estimate for $\widetilde a$ and the $\dot B_{2,1}^{\frac{d}{2}-1}$-regularity estimate for $\widetilde u$. To this end, we introduce the following functional:
\begin{align}\label{G4.2}
\widetilde {\mathcal{X}}(t):=  \|\widetilde  a\|_{\widetilde L_t^\infty( \dot B_{2,1}^{\frac{d}{2}})}+\|\widetilde  u\|_{\widetilde L^\infty_t( \dot B_{2,1}^{\frac{d}{2}-1})}+\|\widetilde u\|_{L^1_t( \dot B_{2,1}^{\frac{d}{2}+1})}.   
\end{align}
We now proceed to the detailed calculations.

\begin{lem}\label{L4.1}
It holds that
\begin{align} \label{G4.3}
\|\widetilde  u\|_{\widetilde L_t^\infty( \dot B_{2,1}^{\frac{d}{2}-1})}+\|\widetilde  u\|_{L_t^1( \dot B_{2,1}^{\frac{d}{2}+1})}\lesssim \|\widetilde a_0\|_{\dot B_{2,1}^\frac{d}{2}}+\|\widetilde u_0\|_{\dot B_{2,1}^{\frac{d}{2}-1}}+(\delta_0+\delta_1)    \widetilde {\mathcal{X}}(t), 
\end{align}
for all $t>0$.
\end{lem}
\begin{proof}
For the  parabolic equation \eqref{G4.1}$_2$, by using Lemma \ref{LA.6}, we derive that
\begin{align}\label{G4.4}
\|\widetilde  u\|_{\widetilde L_t^\infty( \dot B_{2,1}^{\frac{d}{2}-1})}+\|\widetilde  u\|_{L_t^1( \dot B_{2,1}^{\frac{d}{2}+1})}\lesssim&\, \|\widetilde u_0\|_{\dot B_{2,1}^{\frac{d}{2}-1}}+\|\widetilde u\cdot\nabla\bar u\|_{L_t^1( \dot B_{2,1}^{\frac{d}{2}-1})}+\|  u\cdot\nabla\widetilde u\|_{L_t^1( \dot B_{2,1}^{\frac{d}{2}-1})} \nonumber\\
&+\|\widetilde F_1\|_{L_t^1( \dot B_{2,1}^{\frac{d}{2}-1})}+\|\widetilde F_2\|_{L_t^1( \dot B_{2,1}^{\frac{d}{2}-1})}.   
\end{align}
By Lemmas \ref{LA.2}--\ref{LA.4}, the nonlinear terms on the right hand side of \eqref{G4.4} can be estimated as follows:
\begin{align}\label{G4.5}
  \|\widetilde u\cdot\nabla\bar u\|_{L_t^1( \dot B_{2,1}^{\frac{d}{2}-1})}\lesssim&\, \|\widetilde u\|_{L_t^\infty(\dot B_{2,1}^{\frac{d}{2}-1})}\|\bar u\|_{L_t^1(\dot B_{2,1}^\frac{d}{2})}\lesssim \delta_1 \widetilde{\mathcal{X}}(t),\\ \label{G4.6}
  \| u\cdot\nabla\widetilde u\|_{L_t^1( \dot B_{2,1}^{\frac{d}{2}-1})}\lesssim&\, \| u\|_{L_t^\infty(\dot B_{2,1}^{\frac{d}{2}-1})}\|\widetilde u\|_{L_t^1(\dot B_{2,1}^\frac{d}{2})}\lesssim \delta_0 \widetilde{\mathcal{X}}(t),\\ \label{G4.7}
  \|\widetilde F_1\|_{L_t^1(\dot B_{2,1}^{\frac{d}{2}-1})}\lesssim&\, \|f(a )-f(\bar a)\|_{\widetilde L_t^\infty(\dot B_{2,1}^{\frac{d}{2} })}\|\Delta u\|_{L_t^1(\dot B_{2,1}^{\frac{d}{2}-1})}+\|f(\bar a)\|_{\widetilde L_t^\infty(\dot B_{2,1}^{\frac{d}{2} })}\|\Delta\widetilde  u \|_{L_t^1(\dot B_{2,1}^{\frac{d}{2}-1})}
   \nonumber\\
  \lesssim&\,  \|\widetilde a\|_{L_t^\infty(\dot B_{2,1}^\frac{d}{2})}\|  u \|_{L_t^1(\dot B_{2,1}^{\frac{d}{2}+1})}+\|\bar a \|_{L_t^\infty(\dot B_{2,1}^\frac{d}{2})}\| \widetilde u \|_{L_t^1(\dot B_{2,1}^{\frac{d}{2}+1})}\nonumber\\
 \lesssim&\, (\delta_0+\delta_1)\widetilde{\mathcal{X}}(t),\\ \label{G4.8}
    \|\widetilde F_2\|_{L_t^1(\dot B_{2,1}^{\frac{d}{2}-1})} 
  \lesssim&\,  \|\widetilde a\|_{L_t^\infty(\dot B_{2,1}^\frac{d}{2})}\|  u \|_{L_t^1(\dot B_{2,1}^{\frac{d}{2}+1})}+\|\bar a \|_{L_t^\infty(\dot B_{2,1}^\frac{d}{2})}\| \widetilde u \|_{L_t^1(\dot B_{2,1}^{\frac{d}{2}+1})}\nonumber\\
  \lesssim&\, (\delta_0+\delta_1)\widetilde{\mathcal{X}}(t),
\end{align}
where we have   utilized the result from \eqref{TB.3}, which states that 
\begin{align*}
\sup_{(t,x)\in \mathbb R^+ \times\mathbb R^d} |(a,\bar a)(t,x)|\leq \frac{1}{2} \,\Rightarrow \,\frac{1}{2}\leq 1+a\leq \frac{3}{2}\quad \text{and}\quad\frac{1}{2}\leq 1+\bar a\leq \frac{3}{2} .   
\end{align*}
Plugging the estimates \eqref{G4.5}--\eqref{G4.8} into \eqref{G4.4} gives rise to \eqref{G4.3}.
\end{proof}

\begin{lem}\label{L4.2}
It holds that
\begin{align} \label{G4.9}
\|\widetilde  a\|_{\widetilde L_t^\infty( \dot B_{2,1}^{\frac{d}{2}})}\lesssim \|\widetilde a_0\|_{\dot B_{2,1}^\frac{d}{2}}+\|\widetilde u_0\|_{\dot B_{2,1}^{\frac{d}{2}-1}}+(\delta_0+\delta_1)    \widetilde {\mathcal{X}}(t), 
\end{align}
for all $t>0$.
\end{lem}
\begin{proof}
Applying the operator $\dot\Delta_k$ to \eqref{G4.1} and performing a standard $L^2$ energy estimate, one has
\begin{align*}
 \frac{{\rm d}}{{\rm d}t}\|\dot\Delta_k \widetilde a\|_{L^2}^2
  \lesssim&\, \big(\|[u\cdot\nabla, \dot\Delta_k]\widetilde a\|_{L^2}+\|[\widetilde u\cdot\nabla, \dot\Delta_k]\bar a\|_{L^2}+\|\dot\Delta_k(\bar a{\rm div}\, \widetilde u)\|_{L^2}+ 2^k\|\dot\Delta_k\widetilde u\|_{L^2}\big)  \|\dot\Delta_k \widetilde a\|_{L^2}\nonumber\\
&+\big( \|\dot\Delta_k(\widetilde a{\rm div}\,   u)\|_{L^2}+\|{\rm div}\, u\|_{L^\infty}\|\dot\Delta_k \widetilde a\|_{L^2} +\|{\rm div}\, \widetilde u\|_{L^\infty}\|\dot\Delta_k \bar a\|\big)\|\dot\Delta_k \widetilde a\|_{L^2},
\end{align*}
which implies
\begin{align}\label{G4.10}
\|\widetilde a\|_{\widetilde L_t^\infty(\dot B_{2,1}^\frac{d}{2})}\lesssim&\, \|\widetilde a_0\|_{\dot B_{2,1}^\frac{d}{2}}+    \sum_{k\in \mathbb Z}2^{\frac{kd}{2}}\|[u\cdot\nabla,\dot\Delta_k] \widetilde a\|_{L_t^1(L^2)}+    \sum_{k\in \mathbb Z}2^{\frac{kd}{2}}\|[\widetilde u\cdot\nabla,\dot\Delta_k] \bar a\|_{L_t^1(L^2)}\nonumber\\
&+\|\bar a{\rm div}\, \widetilde u\|_{L_t^1(\dot B_{2,1}^\frac{d}{2})}+\|\widetilde a{\rm div}\,   u\|_{L_t^1(\dot B_{2,1}^\frac{d}{2})}+\|\widetilde u\|_{L_t^1(\dot B_{2,1}^{\frac{d}{2}+1})}\nonumber\\
&+\|{\rm div}\, u\|_{L_t^1( L^\infty)}\|\widetilde a\|_{\widetilde L_t^\infty(\dot B_{2,1}^\frac{d}{2})}+\|{\rm div}\,  \widetilde u\|_{L_t^1( L^\infty)}\|\bar a\|_{\widetilde L_t^\infty(\dot B_{2,1}^\frac{d}{2})}.    
\end{align}
It follows from Lemmas \ref{LA.2}--\ref{LA.3} and  \ref{LA.5} that 
\begin{align*}
 \sum_{k\in \mathbb Z}2^{\frac{kd}{2}}\|[u\cdot\nabla,\dot\Delta_k] \widetilde a\|_{L_t^1(L^2)} \lesssim&\,  \|u\|_{L_t^1(\dot B_{2,1}^{\frac{d}{2}+1})}\|\widetilde a\|_{L_t^\infty(\dot B_{2,1}^\frac{d}{2})} \lesssim \delta_0 \widetilde {\mathcal{X}}(t),\\     
 \sum_{k\in \mathbb Z}2^{\frac{kd}{2}}\|[\widetilde u\cdot\nabla,\dot\Delta_k] \bar a\|_{L_t^1(L^2)}\lesssim&\, \|\widetilde u\|_{L_t^1(\dot B_{2,1}^{\frac{d}{2}+1})}\|\bar a\|_{L_t^\infty(\dot B_{2,1}^\frac{d}{2})} \lesssim \delta_1\widetilde {\mathcal{X}}(t),\\     
 \|\bar a{\rm div}\,\widetilde u\|_{L_t^1(\dot B_{2,1}^\frac{d}{2})}\lesssim&\,\|\bar a\|_{L_t^\infty(\dot B_{2,1}^\frac{d}{2})}\|\widetilde u\|_{L_t^1(\dot B_{2,1}^{\frac{d}{2}+1})} \lesssim \delta_1\widetilde {\mathcal{X}}(t),\\ 
  \|\widetilde a{\rm div}\,  u\|_{L_t^1(\dot B_{2,1}^\frac{d}{2})}\lesssim&\,\|\widetilde a\|_{L_t^\infty(\dot B_{2,1}^\frac{d}{2})}\|  u\|_{L_t^1(\dot B_{2,1}^{\frac{d}{2}+1})} \lesssim \delta_0\widetilde {\mathcal{X}}(t),\\ 
  \|{\rm div}\, u\|_{L_t^1( L^\infty)}\|\widetilde a\|_{\widetilde L_t^\infty(\dot B_{2,1}^\frac{d}{2})}\lesssim&\, \|u\|_{L_t^1(\dot B_{2,1}^{\frac{d}{2}+1})} \|\widetilde a\|_{\widetilde L_t^\infty(\dot B_{2,1}^\frac{d}{2})}\lesssim\delta_0\widetilde {\mathcal{X}}(t),\\
  \|{\rm div}\,  \widetilde u\|_{L_t^1( L^\infty)}\|\bar a\|_{\widetilde L_t^\infty(\dot B_{2,1}^\frac{d}{2})}\lesssim&\,\|\widetilde u\|_{L_t^1(\dot B_{2,1}^{\frac{d}{2}+1})} \|\bar a\|_{\widetilde L_t^\infty(\dot B_{2,1}^\frac{d}{2})}\lesssim\delta_1\widetilde {\mathcal{X}}(t).
\end{align*}
Putting all the aforementioned estimates into \eqref{G4.10} yields
\begin{align*} 
\|\widetilde a\|_{\widetilde L_t^\infty(\dot B_{2,1}^\frac{d}{2})}\lesssim&\, \|\widetilde a_0\|_{\dot B_{2,1}^\frac{d}{2}}+ (\delta_0+\delta_1)     \widetilde {\mathcal{X}}(t)+   \|\widetilde u\|_{L_t^1(\dot B_{2,1}^{\frac{d}{2}+1})}.   
\end{align*}
Thanks to the estimate of $\|\widetilde u\|_{L_t^1(\dot B_{2,1}^{\frac{d}{2}+1})}$ established in \eqref{G4.3}, we  further derive \eqref{G4.9}.
\end{proof}
\begin{rem}\label{REM3.1}
In fact, the transport estimate \eqref{A.9} obtained in \cite[Chapter 3]{BCD-Book-2011} does not apply in the present context. Applying the regularity estimate \eqref{A.9} from Lemma \ref{A.7} to equation \eqref{G4.1}$_1$ yields
\begin{align*} 
 \|\widetilde  a\|_{\widetilde L_t^\infty( \dot B_{2,1}^{\frac{d}{2}})}\lesssim   {\rm exp} \Big\{C\|u\|_{L_t^1(\dot B_{2,1}^\frac{d}{2})}   \Big\}     &\, \Big(\|\widetilde  a_0\|_{  \dot B_{2,1}^{\frac{d}{2}}} + \|\bar a{\rm div}\, \widetilde u\|_{L_t^1( \dot B_{2,1}^{\frac{d}{2}})}+ \|\widetilde a{\rm div}\,   u\|_{L_t^1( \dot B_{2,1}^{\frac{d}{2}})} \nonumber\\
  &+ \|\widetilde u\cdot\nabla\bar  a\|_{L_t^1( \dot B_{2,1}^{\frac{d}{2}})}+ {\|\widetilde u\|_{L_t^1( \dot B_{2,1}^{\frac{d}{2}+1})}}    \Big).
\end{align*}
However, the term $\|\widetilde u \cdot \nabla \bar a\|_{L_t^1( \dot B_{2,1}^{\frac{d}{2}})}$ cannot be controlled directly by the $\dot B_{2,1}^\frac{d}{2}$-regularity of $\bar a$ and the $\dot B_{2,1}^{\frac{d}{2}-1}$-regularity of $\widetilde u$.
\end{rem}

With the help of Lemmas \ref{L4.1} and \ref{L4.2}, we proceed to establish the uniform estimate \eqref{error}.
\begin{proof}[Proof of Theorem \ref{T1.3}]
Adding \eqref{G4.3} and \eqref{G4.9} up, we have
 \begin{align*}
  \widetilde {\mathcal{X}}(t)\lesssim   \|\widetilde a_0\|_{\dot B_{2,1}^\frac{d}{2}}+\|\widetilde u_0\|_{\dot B_{2,1}^{\frac{d}{2}-1}}+(\delta_0+\delta_1)    \widetilde {\mathcal{X}}(t), 
 \end{align*}
for all $t>0$. By utilizing the smallness of $\delta_0$ and $\delta_1$, we have for all $t > 0$ that
\begin{align*}
\widetilde{\mathcal{X}}(t) \lesssim \|\widetilde{a}_0\|_{\dot{B}_{2,1}^{\frac{d}{2}}} + \|\widetilde{u}_0\|_{\dot{B}_{2,1}^{\frac{d}{2}-1}}.
\end{align*}
Thus, the proof of Theorem \ref{T1.3} is completed.
\end{proof}

\section{Optimal time decay rate  of  the pressureless Navier-Stokes system}

\subsection{The upper bound estimate for $u$ under the boundedness condition}
In this subsection, we aim to establish the upper bound of the decay estimates for $u$. First, we analyze the propagation  of $\dot B_{2,\infty}^{\sigma_0}$ in the low-frequency regime.
\begin{lem}\label{L5.1}
Let $(a,u)$ be the  global strong solution to the Cauchy problem \eqref{A4} given by Theorem \ref{T1.2}. Then, under the assumptions of
Theorem \ref{T1.4}, we have
\begin{align}\label{G5.1}
\mathcal{X}_{\ell,\sigma_0}(t):= \|u\|^\ell_{\widetilde L_t^\infty(\dot B_{2,\infty}^{\sigma_0})}+  \|u\|^\ell_{  L_t^1(\dot B_{2,\infty}^{\sigma_0+2})} \leq C\delta_{*},
\end{align}
for all $t > 0$, where $\delta_*$ is given in \eqref{TC.3} and $C > 0$ is a constant independent of time.
\end{lem}
\begin{proof}
Multiplying \eqref{G3.6} by $2^{k\sigma_0}$ and taking the supremum on  both $[0,t]$ and   $k \leq 0$, one has 
\begin{align}\label{G5.2}
&\|u\|_{\widetilde L_t^\infty(\dot B_{2,\infty}^{ \sigma_0})}^\ell+\|u\|_{L_t^1(\dot B_{2,\infty}^{ \sigma_0+2})}^\ell\nonumber\\
&\quad\lesssim  \|u_0\|_{\dot B_{2,\infty}^{ \sigma_0}}^\ell +\|u\cdot\nabla u\|_{L_t^1(\dot B_{2,\infty}^{ \sigma_0})}^\ell+ \|f(a)\Delta u\|_{L_t^1(\dot B_{2,\infty}^{ \sigma_0})}^\ell+  \|f(a)\nabla{\rm div}\, u\|_{L_t^1(\dot B_{2,\infty}^{ \sigma_0})}^\ell.    
\end{align}

By applying an argument analogous to that in Lemma \ref{L3.2} and building upon Lemmas \ref{LA.2}--\ref{LA.4}, we conclude that  
\begin{align}\label{G5.3}
 \|u\cdot\nabla u\|_{L_t^1(\dot B_{2,\infty}^{ \sigma_0})}^\ell\lesssim&\, \|u\|_{ L_t^2(\dot B_{2,1}^{\frac{d}{2}})} \|u\|_{L_t^2(\dot B_{2,\infty}^{\sigma_0+1})}\nonumber\\
 \lesssim&\,   \|u\|_{ L_t^1(\dot B_{2,1}^{\frac{d}{2}+1})}^{\frac{1}{2}} \|u\|_{ \widetilde L_t^\infty(\dot B_{2,1}^{\frac{d}{2}-1})}^{\frac{1}{2}} \|u\|_{ L_t^1(\dot B_{2,\infty}^{ \sigma_0+2})}^{\frac{1}{2}} \|u\|_{ \widetilde L_t^\infty(\dot B_{2,\infty}^{\sigma_0})}^{\frac{1}{2}} \nonumber\\
 \lesssim&\, \mathcal{X}(t)\Big(\|u^\ell\|_{ L_t^1(\dot B_{2,\infty}^{ \sigma_0+2})} +\|u^h\|_{ L_t^1(\dot B_{2,\infty}^{ \sigma_0+2})}+  \|u^\ell\|_{ \widetilde L_t^\infty(\dot B_{2,\infty}^{\sigma_0})}  +  \|u^h\|_{ \widetilde L_t^\infty(\dot B_{2,\infty}^{\sigma_0})}       \Big)\nonumber\\
 \lesssim&\, \mathcal{X}(t)\big( \mathcal{X}_{\ell,\sigma_0}(t)+  \mathcal{X}(t)     \big),
\end{align}
and 
\begin{align}\label{G5.4}
\|f(a)\Delta u\|_{L_t^1(\dot B_{2,\infty}^{ \sigma_0})}+  \|f(a)\nabla{\rm div}\, u\|_{L_t^1(\dot B_{2,\infty}^{ \sigma_0})} 
\lesssim&\, \|f(a)\|_{\widetilde L_t^\infty(\dot B_{2,1}^\frac{d}{2})}\|u\|_{L_t^1(\dot B_{2,\infty}^{\sigma_0+2})}\nonumber\\
\lesssim&\, \|a\|_{\widetilde L_t^\infty(\dot B_{2,1}^\frac{d}{2})}\Big( \|u^\ell\|_{L_t^1(\dot B_{2,\infty}^{\sigma_0+2})}+\|u^h\|_{L_t^1(\dot B_{2,\infty}^{\sigma_0+2})} \Big)\nonumber\\
\lesssim&\,    \mathcal{X}(t)\big( \mathcal{X}_{\ell,\sigma_0}(t)+  \mathcal{X}(t)     \big).
\end{align}
Inserting the estimates \eqref{G5.3} and \eqref{G5.4} into \eqref{G5.2} gives
\begin{align*}
 \mathcal{X}_{\ell,\sigma_0}(t)\lesssim   \|u_0\|_{\dot B_{2,\infty}^{ \sigma_0}}^\ell +  \mathcal{X}(t)\big( \mathcal{X}_{\ell,\sigma_0}(t)+  \mathcal{X}(t)  \big), 
\end{align*}
which, together with the smallness of $\mathcal{X}(t)$ and  the estimate \eqref{TB.3}, implies  that
\begin{align*}
 \mathcal{X}_{\ell,\sigma_0}(t)\lesssim   \|u_0\|_{\dot B_{2,\infty}^{ \sigma_0}}^\ell + \|a_0\|_{\dot B_{2,1}^\frac{d}{2}}+\|u_0\|_{\dot B_{2,1}^{\frac{d}{2}-1}}\lesssim  \delta_*.
\end{align*}
Thus, we derive \eqref{G5.1}, thereby completing the proof of Lemma \ref{L5.1}.
\end{proof}
Motivated by \cite{LS-SIMA-2023}, we now establish the time-weighted estimates for both the low-frequency and high-frequency components of $u$. For $M > \max\left\{\frac{1}{2}\big(\frac{d}{2} + 1 - \sigma_0\big), 1\right\}$, we introduce a time-weighted functional $\mathcal{X}_{M}(t)$ defined by:
\begin{align}\label{G5.5}
 \mathcal{X}_{M}(t):=\|\tau^M u\|_{\widetilde L_t^\infty(\dot B_{2,1}^{{\frac{d}{2}}-1})} +\|\tau^M u\|_{L_t^1(\dot B_{2,1}^{\frac{d}{2}+1})}.   
\end{align}

\begin{lem}\label{L5.2}
Let $(a,u)$ be the  global strong solution to the Cauchy problem \eqref{A4} given by Theorem \ref{T1.2}. Then, under the assumptions of
Theorem \ref{T1.4},  for any $t>0$ and  $M > \max\left\{\frac{1}{2}\big(\frac{d}{2} + 1 - \sigma_0\big), 1\right\}$, it holds that
\begin{align}\label{G5.6}
 \|\tau^M u\|_{\widetilde L_t^\infty(\dot B_{2,1}^{{\frac{d}{2}}-1})}^\ell+\|\tau^M u\|_{L_t^1(\dot B_{2,1}^{\frac{d}{2}+1})}^\ell\lesssim&   \big(\eps+\mathcal{X}(t)\big) \mathcal{X}_{M}(t)+\frac{\mathcal{X}(t)+\mathcal{X}_{\ell,\sigma_0}(t)}{\eps} t^{M-\frac{1}{2}(\frac{d}{2}-1-\sigma_0)},  
\end{align}
where $\eps>0$ is a constant to be determined later. Here, $\mathcal{X}(t)$, $\mathcal{X}_{\ell,\sigma_0}(t)$, and $\mathcal{X}_M(t)$ are defined through \eqref{defnx}, \eqref{G5.1}, and \eqref{G5.5}, respectively.
\end{lem}
\begin{proof}
Multiplying the inequality \eqref{G3.5} by $t^{ M}$, one gets 
\begin{align}\label{G5.7}
&\frac{{\rm d}}{{\rm d}t}(t^M \|\dot\Delta_k u\|_{L^2}^2)+\lambda_1t^M2^{2k}\|\dot\Delta_k u\|_{L^2}^2- {Mt^{M-1}\|\dot\Delta_k u\|_{L^2}^2} \nonumber\\
&\quad\lesssim t^M \|\dot\Delta_k u\|_{L^2}\big(\|\dot\Delta_k(u\cdot\nabla u)\|_{L^2}+\|\dot\Delta_k(f(a)\Delta u)\|_{L^2}+\|\dot\Delta_k(f(a)\nabla{\rm div}\, u)\|_{L^2}        \big),
\end{align}
for some constant $\lambda_1>0$.  Integrating \eqref{G5.7} over $[0,t]$ and taking the square root of both sides of the resulting inequality, we have  that for any $k\leq 0$,
\begin{align}\label{G5.8}
&t^M \|\dot\Delta_k u\|_{L^2}+2^{2k}\int_0^t \tau^M \|\dot\Delta_k u\|_{L^2}{\rm d}\tau   \nonumber\\
& \quad\lesssim\int_0^t\tau^M\big(\|\dot\Delta_k(u\cdot\nabla u)\|_{L^2}+\|\dot\Delta_k(f(a)\Delta u)\|_{L^2}+\|\dot\Delta_k(f(a)\nabla{\rm div}\, u)\|_{L^2}        \big){\rm d}\tau\nonumber\\
& \quad\quad+ \int_0^t \tau^{M-1} \|\dot\Delta_k u\|_{L^2}{\rm d}\tau.
\end{align}
By multiplying \eqref{G5.8} by $2^{k(\frac{d}{2}-1)}$ and taking the summation over $k\leq 0$, we deduce that
\begin{align}\label{G5.9}
\|\tau^M u\|_{\widetilde L_t^\infty(\dot B_{2,1}^{\frac{d}{2}-1})}^\ell+\|\tau^M u\|_{  L_t^1(\dot B_{2,1}^{\frac{d}{2}+1})}^\ell    \lesssim&\,\|\tau^Mf(a)\Delta u\|_{L_t^1(\dot B_{2,1}^{\frac{d}{2}-1})}^\ell+\|\tau^Mf(a)\nabla{\rm div}\,u\|_{L_t^1(\dot B_{2,1}^{\frac{d}{2}-1})}^\ell\nonumber\\
&+\|\tau^M u\cdot \nabla u\|_{L_t^1(\dot B_{2,1}^{\frac{d}{2}-1})}^\ell+  \int_0^t \tau^{M-1}\|u\|_{\dot B_{2,1}^{\frac{d}{2}-1}}^\ell {\rm d}\tau.  
\end{align}
Making use of an argument analogous to that in Lemma \ref{L3.2} and performing direct calculations, we get
\begin{align}\label{G5.10}
  \|\tau^Mf(a)\Delta u\|_{L_t^1(\dot B_{2,1}^{\frac{d}{2}-1})} +\|\tau^Mf(a)\nabla{\rm div}\,u\|_{L_t^1(\dot B_{2,1}^{\frac{d}{2}-1})}  \lesssim&\, \|a\|_{\widetilde L_t^\infty(\dot B_{2,1}^\frac{d}{2})}\| \tau^Mu\|_{L_t^1(\dot B_{2,1}^{\frac{d}{2}+1})}\nonumber\\
  \lesssim&\, \mathcal{X}(t)\mathcal{X}_{M}(t),\\\label{G5.11}
  \|\tau^M u\cdot \nabla u\|_{L_t^1(\dot B_{2,1}^{\frac{d}{2}-1})}\lesssim&\,  \|u\|_{\widetilde L_t^\infty(\dot B_{2,1}^{\frac{d}{2}-1})}\| \tau^Mu\|_{L_t^1(\dot B_{2,1}^{\frac{d}{2}+1})}\nonumber\\
   \lesssim&\, \mathcal{X}(t)\mathcal{X}_{M}(t).
\end{align}
For the remaining term on the right-hand side of \eqref{G5.9}, we decompose it into two parts:  
\begin{align}\label{G5.12}
\int_0^t \tau^{M-1}\|u\|_{\dot B_{2,1}^{\frac{d}{2}-1}}  {\rm d}\tau\lesssim   \int_0^t \tau^{M-1}\|u^\ell\|_{\dot B_{2,1}^{\frac{d}{2}-1}}  {\rm d}\tau  +\int_0^t \tau^{M-1}\|u^h\|_{\dot B_{2,1}^{\frac{d}{2}-1}}  {\rm d}\tau\equiv: I_1+I_2.
\end{align}
On the one hand, by leveraging \eqref{A.1} in Lemma \ref{LA.2},
we arrive at
\begin{align}\label{G5.13}
I_1\lesssim \int_0^t \tau^{M-1} \|u^\ell\|_{\dot B_{2,\infty}^{\sigma_0}}^{1-\theta_1} \|u^\ell\|_{\dot B_{2,1}^{\frac{d}{2}+1}}^{\theta_1} {\rm d}\tau \lesssim&\, \bigg(\int_0^t \tau^{M-\frac{1}{1-\theta_1}}{\rm d}\tau  \bigg)^{1-\theta_1} \|u^\ell\|_{\widetilde L_t^\infty(\dot B_{2,\infty}^{\sigma_0})}^{1-\theta_1} \|\tau^M u^\ell\|_{L_t^1(\dot B_{2,1}^{\frac{d}{2}+1})}^{\theta_1}\nonumber\\
\lesssim&\,\Big( t^{M-\frac{1}{2}(\frac{d}{2}-1-\sigma_0)} \|u^\ell\|_{\widetilde L_t^\infty(\dot B_{2,\infty}^{\sigma_0})}\Big)^{1-\theta_1} \Big(\|\tau^Mu\|^\ell_{L_t^1(\dot B_{2,1}^{\frac{d}{2}+1})}\Big)^{\theta_1}\nonumber\\
\lesssim&\, \eps {\mathcal{X}_M(t)}+\frac{\mathcal{X}_{\ell,\sigma_0}}{\eps} t^{M-\frac{1}{2}(\frac{d}{2}-1-\sigma_0)},
\end{align}
where the constant $\theta_1=\frac{\frac{d}{2}-1-\sigma_0}{\frac{d}{2}+1-\sigma_0}\in (0,1)$.
On the other hand, it is easy to find that
\begin{align}\label{G5.14}
I_2\lesssim&\,    \Big( t^{M-\frac{1}{2}(\frac{d}{2}-1-\sigma_0)} \|u^h\|_{\widetilde L_t^\infty(\dot B_{2,1}^{\frac{d}{2}-1})}\Big)^{1-\theta_1} \Big(\|\tau^Mu\|^h_{L_t^1(\dot B_{2,1}^{\frac{d}{2}-1})}\Big)^{\theta_1}\nonumber\\
\lesssim&\, \eps {\mathcal{X}_M(t)}+\frac{\mathcal{X}(t)}{\eps} t^{M-\frac{1}{2}(\frac{d}{2}-1-\sigma_0)}. 
\end{align}
Substituting the estimates \eqref{G5.13} and \eqref{G5.14} into \eqref{G5.12} yields
\begin{align}\label{G5.15}
\int_0^t \tau^{M-1}\|u\|_{\dot B_{2,1}^{\frac{d}{2}-1}}  {\rm d}\tau\lesssim   \eps {\mathcal{X}_M(t)}+\frac{\mathcal{X}(t)+\mathcal{X}_{\ell,\sigma_0}(t)}{\eps} t^{M-\frac{1}{2}(\frac{d}{2}-1-\sigma_0)} .
\end{align}
Putting the estimates \eqref{G5.10}, \eqref{G5.11} and \eqref{G5.15} in \eqref{G5.9}, we consequently obtain \eqref{G5.6}.    
\end{proof}

\begin{lem}\label{L5.3}
Let $(a,u)$ be the  global strong solution to the Cauchy problem \eqref{A4} given by Theorem \ref{T1.2}. Then, under the assumptions of
Theorem \ref{T1.4},  for any $t>0$ and  $M > \max\left\{\frac{1}{2}\big(\frac{d}{2} + 1 - \sigma_0\big), 1\right\}$, it holds that
\begin{align} \label{G5.16}
 \|\tau^M u\|_{\widetilde L_t^\infty(\dot B_{2,1}^{{\frac{d}{2}}-1})}^h+\|\tau^M u\|_{L_t^1(\dot B_{2,1}^{\frac{d}{2}+1})}^h\lesssim  \big(\eps+\mathcal{X}(t)\big) \mathcal{X}_{M}(t)+\frac{\mathcal{X}(t)}{\eps} t^{M-\frac{1}{2}(\frac{d}{2}-1-\sigma_0)},  
\end{align}
where $\eps>0$ is a constant to be determined later. Here, $\mathcal{X}(t)$, $\mathcal{X}_{\ell,\sigma_0}(t)$, and $\mathcal{X}_M(t)$ are defined via \eqref{defnx}, \eqref{G5.1}, and \eqref{G5.5}, respectively.
\end{lem}
\begin{proof}
Multiplying the inequality \eqref{G3.11} by $t^M$, we obtain that for any $k\geq -1$,
\begin{align}\label{G5.17}
&\frac{{\rm d}}{{\rm d}t}(t^M \|\dot\Delta_k u\|_{L^2}^2)+\lambda_2 t^M2^{2k}\|\dot\Delta_k u\|_{L^2}^2- {Mt^{M-1}\|\dot\Delta_k u\|_{L^2}^2} \nonumber\\
&\quad\lesssim t^M \|\dot\Delta_k u\|_{L^2}\big(\|[u\cdot\nabla,\dot\Delta_k] u\|_{L^2}+\|\dot\Delta_k(f(a)\Delta u)\|_{L^2}+\|\dot\Delta_k(f(a)\nabla{\rm div}\, u)\|_{L^2}        \big) \nonumber\\
&\quad\quad+t^M\|{\rm div}\, u\|_{L^\infty}\|\dot\Delta_k u\|_{L^2}^2,
\end{align}
for some constant $\lambda_2>0$. 
Integrating the inequality \eqref{G5.17} over the interval $[0,t]$, taking the square root of both sides of the resulting inequality, multiplying \eqref{G5.8} by $2^{k(\frac{d}{2}-1)}$ and taking the summation over $k\geq 1$, we   deduce that   
\begin{align}\label{G5.18}
&\|\tau^M u\|_{\widetilde L_t^\infty(\dot B_{2,1}^{\frac{d}{2}-1})}^h+\|\tau^Mu\|_{L_t^1(\dot B_{2,1}^{\frac{d}{2}+1})}^h\nonumber\\
&\quad\lesssim   \|{\rm div}\,u\|_{L_t^1(L^\infty)}\|\tau^M u \|_{\widetilde L_t^\infty(\dot B_{2,1}^{\frac{d}{2}-1})}^h+\sum_{k\geq -1}2^{k(\frac{d}{2}-1)}\|[u\cdot\nabla,\dot\Delta_k] \tau^Mu\|_{L_t^1(L^2)} \nonumber\\
&\quad\quad+ \|\tau^Mf(a)\Delta u\|_{L_t^1(\dot B_{2,1}^{\frac{d}{2}-1})}^h+  \|\tau^Mf(a)\nabla{\rm div}\, u\|_{L_t^1(\dot B_{2,1}^{\frac{d}{2}-1})}^h+\int_0^t\tau^{M-1} \|u\|_{\dot B_{2,1}^{\frac{d}{2}-1}}^h{\rm d}\tau.      
\end{align}
It follows from Lemmas \ref{LA.2}--\ref{LA.5} that 
\begin{align}\label{G5.19}
\|{\rm div}\,u\|_{L_t^1(L^\infty)}\|\tau^M u \|_{\widetilde L_t^\infty(\dot B_{2,1}^{\frac{d}{2}-1})}^h\lesssim&\, \|u\|_{L_t^1(\dot B_{2,1}^{\frac{d}{2}+1})} \|\tau^M u\|_{\widetilde L_t^\infty(\dot B_{2,1}^{\frac{d}{2}-1})} \lesssim\mathcal{X}(t)\mathcal{X}_{M}(t),  \\ \label{G5.20}
\sum_{k\geq -1}2^{k(\frac{d}{2}-1)}\|[u\cdot\nabla,\dot\Delta_k] \tau^Mu\|_{L_t^1(L^2)}\lesssim&\, \|u\|_{\widetilde L^\infty_t(\dot B_{2,1}^{\frac{d}{2}-1})}\|\tau^Mu\|_{L_t^1(\dot B_{2,1}^{\frac{d}{2}+1})}\lesssim\mathcal{X}(t)\mathcal{X}_{M}(t).
\end{align} 
Similar to \eqref{G5.12} and \eqref{G5.13}, we deduce that
\begin{align}\label{G5.21}
\int_0^t\tau^{M-1} \|u\|_{\dot B_{2,1}^{\frac{d}{2}-1}}^h{\rm d}\tau\lesssim &\,  \Big( t^{M-\frac{1}{2}(\frac{d}{2}-1-\sigma_0)} \|u\|^h_{\widetilde L_t^\infty(\dot B_{2,1}^{\frac{d}{2}-1})}\Big)^{1-\theta_2} \Big(\|\tau^Mu\|^h_{L_t^1(\dot B_{2,1}^{\frac{d}{2}+1})}\Big)^{\theta_2}\nonumber\\
\lesssim&\, \eps  {\mathcal{X}_M(t)}+\frac{\mathcal{X}(t)}{\eps} t^{M-\frac{1}{2}(\frac{d}{2}-1-\sigma_0)},  
\end{align}
where the constant $\theta_2=\frac{\frac{d}{2}-1-\sigma_0}{\frac{d}{2}+1-\sigma_0}\in (0,1)$.
Plugging the estimates \eqref{G5.10} and \eqref{G5.19}--\eqref{G5.21} into \eqref{G5.18}, we directly obtain \eqref{G5.16}.
\end{proof}
Next, we proceed to prove Theorem \ref{T1.4}.
\begin{proof}[Proof of Theorem \ref{T1.4}]
By combining \eqref{G5.6} and \eqref{G5.16} with \eqref{G5.5}, one has 
\begin{align}\label{G5.22}
\mathcal{X}_{M}(t)\lesssim  \big(\eps+\mathcal{X}(t)\big) \mathcal{X}_{M}(t)+\frac{\mathcal{X}(t)+\mathcal{X}_{\ell,\sigma_0}(t)}{\eps} t^{M-\frac{1}{2}(\frac{d}{2}-1-\sigma_0)},    
\end{align}
for any $t>0$ and  $M > \max\left\{\frac{1}{2}\big(\frac{d}{2} +1 - \sigma_0\big), 1\right\}$. By choosing a sufficiently small constant $\eps>0$ in \eqref{G5.22} and exploiting the smallness of $\mathcal{X}(t)$ together with \eqref{G5.1}, we derive that
\begin{align}\label{G5.23}
 \mathcal{X}_{M}(t)\lesssim   \delta_{*} t^{M-\frac{1}{2}(\frac{d}{2}-1-\sigma_0)},
\end{align}
for any $t>0$, which implies that
\begin{align}\label{G5.24}
\|u(t)\|_{\dot B_{2,1}^{\frac{d}{2}-1}}\lesssim\delta_{*} t^{-\frac{1}{2}(\frac{d}{2}-1-\sigma_0)},    
\end{align}
for any $t\geq 1$.    Since $\|u(t)\|_{\dot B_{2,1}^{\frac{d}{2}-1}}$ is uniformly bounded by $\delta_{*}$, it follows from \eqref{G5.24} that 
\begin{align}\label{G5.25}
\|u(t)\|_{\dot B_{2,1}^{\frac{d}{2}-1}}\lesssim\delta_{*} (1+ t)^{-\frac{1}{2}(\frac{d}{2}-1-\sigma_0)},    
\end{align}
for any $t\geq 1$, which  together with \eqref{A.1} in Lemma \ref{LA.2} yields  
\begin{align*}
\|u(t)\|_{\dot B_{2,1}^{\sigma}}\lesssim&\,  \|u^\ell(t)\|_{\dot B_{2,1}^{\sigma}} +   \|u^h(t)\|_{\dot B_{2,1}^{\sigma}}\nonumber\\
\lesssim&\,  \|u^\ell(t)\|_{\dot B_{2,\infty}^{\sigma_0}}^{  \frac{{\frac{d}{2}-1-\sigma}}{\frac{d}{2}-1-\sigma_0} }  \|u^\ell(t)\|_{\dot B_{2,1}^{\frac{d}{2}-1}}^{  \frac{\sigma-\sigma_0}{\frac{d}{2}-1-\sigma_0} }+\|u^h(t)\|_{\dot B_{2,1}^{\frac{d}{2}-1}}\nonumber\\
\lesssim&\, \delta_*^{  \frac{{\frac{d}{2}-1-\sigma}}{\frac{d}{2}-1-\sigma_0} }\Big(\delta_*\|u^\ell(t)\|_{\dot B_{2,\infty}^{\frac{d}{2}-1}}\Big)^{ \frac{\sigma-\sigma_0}{\frac{d}{2}-1-\sigma_0}}+\delta_{*} (1+t)^{-\frac{1}{2}(\frac{d}{2}-1-\sigma_0)}\nonumber\\
\lesssim&\, \delta_* (1+t)^{-\frac{1}{2}(\sigma-\sigma_0)},
\end{align*}
for any $\sigma\in(\sigma_0,\frac{d}{2}-1]$. Then \eqref{TC.2} follows, and it remains to prove \eqref{TC.4}.
From \eqref{TB.3}, it follows that there exists a time $t_1>0$ such that
\begin{align*}
\|u(t_1)\|_{\dot B_{2,1}^{\frac{d}{2}+1}} \lesssim  \|a_0\|_{\dot B_{2,1}^\frac{d}{2}}+\|u_0\|_{\dot B_{2,1}^{\frac{d}{2}-1}}.  
\end{align*}
For simplicity of notation, we set $t_1 = 1$. 
Multiplying \eqref{A4}$_2$ by $t^M$ gives
\begin{align}\label{GNJK5.26}
&\partial_t(t^M u)-\mu \Delta(t^{M} u)-(\mu+\nu)\nabla{\rm div}\,(t^Mu)\nonumber\\
&\quad= Mt^{M-1} u-t^M u\cdot\nabla u+t^Mf(a)\mu \Delta u+t^M f(a)(\mu+\nu){\nabla {\rm div}\, u}.  
\end{align}
By combining \eqref{GNJK5.26} with the maximal regularity estimate for $u$ established in Lemma \ref{LA.6} at high frequencies over the interval $[1, t]$, one has  
\begin{align*}
&\|\tau^M u(t)\|_{\widetilde L^\infty(1,t;\dot B_{2,1}^{\frac{d}{2}+1})}^h\nonumber\\
\lesssim&\, \|u(1)\|^h_{\dot B_{2,1}^{\frac{d}{2}+1}} +\|\tau^{M-1 } u\|_{\widetilde L^\infty(1,t;\dot B_{2,1}^{\frac{d}{2}-1})}^h   +\|\tau^{M } u\cdot\nabla u\|_{\widetilde L^\infty(1,t;\dot B_{2,1}^{\frac{d}{2}-1})}^h \nonumber\\
& +\|\tau^{M} f(a)\Delta u\|_{\widetilde L^\infty(1,t;\dot B_{2,1}^{\frac{d}{2}-1})}^h +\|\tau^{M} f(a)\nabla{\rm div}\,u\|_{\widetilde L^\infty(1,t;\dot B_{2,1}^{\frac{d}{2}-1})}^h\nonumber\\
\lesssim&\, \|u(1)\|^h_{\dot B_{2,1}^{\frac{d}{2}+1}} +\mathcal{X}_{M}(t)+ \|u\|_{\widetilde L^\infty(1,t;\dot B_{2,1}^{\frac{d}{2}-1})}\Big ( \|\tau^M u\|_{\widetilde L^\infty(1,t;\dot B_{2,1}^{\frac{d}{2}+1})}^\ell+  \|\tau^M u\|_{\widetilde L^\infty(1,t;\dot B_{2,1}^{\frac{d}{2}+1})}^h  \Big) \nonumber\\
& +\|a\|_{\widetilde L^\infty(1,t;\dot B_{2,1}^{\frac{d}{2} })}\Big ( \|\tau^M u\|_{\widetilde L^\infty(1,t;\dot B_{2,1}^{\frac{d}{2}+1})}^\ell+  \|\tau^M u\|_{\widetilde L^\infty(1,t;\dot B_{2,1}^{\frac{d}{2}+1})}^h  \Big)\nonumber\\
\lesssim&\,\delta_0\|\tau^M u(t)\|_{\widetilde L^\infty(1,t;\dot B_{2,1}^{\frac{d}{2}+1})}^h+\delta_* t^{M-\frac{1}{2}(\frac{d}{2}-1-\sigma_0)},
\end{align*}
which together with the smallness of $\delta_0$ gives
\begin{align*}
\|\tau^M u(t)\|_{\widetilde L^\infty(1,t;\dot B_{2,1}^{\frac{d}{2}+1})}^h\lesssim \delta_* t^{M-\frac{1}{2}(\frac{d}{2}-1-\sigma_0)},    
\end{align*}
for any $t\geq 1$. Consequently, it holds that
\begin{align*}
\| u(t)\|_{ \dot B_{2,1}^{\frac{d}{2}+1} }^h\lesssim \delta_* t^{-\frac{1}{2}(\frac{d}{2}-1-\sigma_0)}, 
\end{align*}
for any $t\geq 1$. 
Hence, we complete the proof of Theorem \ref{T1.4}.
\end{proof}

\subsection{The upper bound estimate for $u$ under the smallness condition}
In this subsection, we establish the decay estimates for higher-order spatial derivatives of $u$ under the smallness condition of $u_0^\ell$ in $\dot B_{2,\infty}^{\sigma_0}$, where $\sigma_0\in[-\frac{d}{2},\frac{d}{2}-1)$. These estimates yield improved results compared to those in Theorem \ref{T1.4}, as stated in Theorem \ref{T1.1}. The proof is mainly adapted from \cite[Section 5.2]{Danchin-2018}. 
First, we introduce the time-weighted energy functional  $\mathcal{D}(t)$:
\begin{align}\label{NJKG4.27}
\mathcal{D}(t):=\sup_{\sigma\in[\sigma_0+\theta,\frac{d}{2}+1]} \|\langle\tau\rangle^{\frac{1}{2}(\sigma-\sigma_0)} u\|_{L_t^\infty(\dot B_{2,1}^{\sigma})}^\ell+\|\langle\tau\rangle^\alpha u\|_{\widetilde L_t^\infty(\dot B_{2,1}^{\frac{d}{2}-1})}^h+\| \tau ^\alpha u\|_{\widetilde L_t^\infty(\dot B_{2,1}^{\frac{d}{2}+1})}^h,
\end{align}
where  $\langle t\rangle:=\sqrt{1+t^2}$ and $\alpha:=\frac{1}{2}(\frac{d}{2}+1-\sigma_0 )$, and $\theta\in(0,1]$ is a sufficiently small constant. 

Now we provide the estimate of $u$ in both low and high frequencies, as detailed in the following lemmas.
\begin{lem}\label{NJKL4.4}
Let $(a,u)$ be the global strong solution to the Cauchy problem \eqref{A4} given by Theorem \ref{T1.2}. Then, under the assumptions of Theorem \ref{T1.1}, it holds that for any $t>0$,
\begin{align}\label{NJKG4.28}
 \sup_{\sigma\in[\sigma_0+\theta,\frac{d}{2}+1]} \big\|\langle\tau\rangle^{\frac{1}{2}(\sigma-\sigma_0)} u\big\|_{L_t^\infty(\dot B_{2,1}^{\sigma})}^\ell\lesssim  \delta_*+\delta_*\mathcal{D}(t)+\mathcal{D}^2(t),  
\end{align}
where $\delta_{*}$ and $\mathcal{D}(t)$ are defined via                 \eqref{TC.3} and \eqref{NJKG4.27}, respectively.
\end{lem}
\begin{proof}
Using Gronwall's inequality to \eqref{G3.6}, we arrive at
\begin{align}\label{NJKG4.29}
\|\dot\Delta_k u\|_{L^2}\lesssim&\, e^{-2^{2k}t}   \|\dot \Delta_k u_0\|_{L^2}+\int_0^t e^{-2^{2k}(t-\tau)} \|\dot\Delta_k(u\cdot\nabla u)\|_{L^2}{\rm d}\tau\nonumber\\
&+\int_0^t e^{-2^{2k}(t-\tau)} \big(\|\dot\Delta_k(f(a)\Delta u)\|_{L^2}+\|\dot\Delta_k(f(a)\nabla{\rm div}\, u)\|_{L^2}        \big){\rm d}\tau,
\end{align}
for any $k\leq 0$, which implies that
\begin{align}\label{NJKG4.30}
\|u\|_{\dot B_{2,1}^\sigma}^\ell\lesssim&\,\int_0^t \langle t-\tau\rangle^{-\frac{1}{2}(\sigma-\sigma_0)} \Big(\|u\cdot\nabla u\|_{\dot B_{2,\infty}^{\sigma_0}}^\ell+\|f(a)\Delta u\|_{\dot B_{2,\infty}^{\sigma_0}}^\ell+\|f(a)\nabla {\rm div}\, u\|_{\dot B_{2,\infty}^{\sigma_0}}^\ell     \Big) {\rm d}\tau \nonumber\\
&+\langle t\rangle^{-\frac{1}{2}(\sigma-\sigma_0)}\|u_0\|_{\dot B_{2,\infty}^{\sigma_0}}^\ell,   
\end{align}
for any $\sigma>\sigma_0$, where we have used the following inequality (see \cite[Lemma 2.35]{BCD-Book-2011}):
\begin{align*}
\sup_{t>0}\sum_{k\in\mathbb Z}t^s 2^{ks}e^{-ct2^{2k}}\lesssim 1, \quad \forall s>0.     
\end{align*}
To derive the estimate of ${\|u\|_{\dot B_{2,1}^\sigma}^\ell}$, we consider the cases $t \leq 2$ and $t \geq 2$ separately.
For the case $t\leq 2$, we utilize the fact that $\langle t\rangle\backsim1$ to get
\begin{align}\label{NJKG4.31}
\int_0^t \langle t-\tau\rangle^{-\frac{1}{2}(\sigma-\sigma_0)}\|u\cdot\nabla u\|_{\dot B_{2,\infty}^{\sigma_0} }^\ell{\rm d}\tau \lesssim&\,   \|u\|_{ L_t^1(\dot B_{2,1}^{\frac{d}{2}+1})}^{\frac{1}{2}} \|u\|_{ \widetilde L_t^\infty(\dot B_{2,1}^{\frac{d}{2}-1})}^{\frac{1}{2}} \|u\|_{ L_t^1(\dot B_{2,\infty}^{ \sigma_0+2})}^{\frac{1}{2}} \|u\|_{ \widetilde L_t^\infty(\dot B_{2,\infty}^{\sigma_0})}^{\frac{1}{2}} \nonumber\\
 \lesssim&\, \mathcal{X}(t)\big( \mathcal{X}_{\ell,\sigma_0}(t)+  \mathcal{X}(t)     \big) \langle t\rangle^{-\frac{1}{2}(\sigma-\sigma_0)} \nonumber\\   
 \lesssim&\,\delta_{*}\langle t\rangle^{-\frac{1}{2}(\sigma-\sigma_0)} .
\end{align}
For the case $t \geq 2$, we split the integration into two parts:
\begin{align*} 
&\int_0^t \langle t-\tau\rangle^{-\frac{1}{2}(\sigma-\sigma_0)}\|u\cdot\nabla u\|_{\dot B_{2,\infty}^{\sigma_0} }^\ell{\rm d}\tau\nonumber\\
&=\,\int_0^1  \langle t-\tau\rangle^{-\frac{1}{2}(\sigma-\sigma_0)}\|u\cdot\nabla u\|_{\dot B_{2,\infty}^{\sigma_0} }^\ell{\rm d}\tau+  \int_1^t  \langle t-\tau\rangle^{-\frac{1}{2}(\sigma-\sigma_0)}\|u\cdot\nabla u\|_{\dot B_{2,\infty}^{\sigma_0} }^\ell{\rm d}\tau=:J_1+J_2.
\end{align*}
On the one hand, when $\tau\in[0,1]$, it holds that $\langle t-\tau\rangle \sim \langle t\rangle$. Consequently, we derive that
\begin{align*} 
J_1\lesssim \langle t\rangle^{-\frac{1}{2}(\sigma-\sigma_0)} \int_0^1 \|u\|_{\dot B_{2,1}^{\frac{d}{2}}}\|u\|_{\dot B_{2,\infty}^{\sigma_0+1}} {\rm d}\tau  \lesssim  \mathcal{X}(t)\big( \mathcal{X}_{\ell,\sigma_0}(t)+  \mathcal{X}(t)     \big) \langle t\rangle^{-\frac{1}{2}(\sigma-\sigma_0)}\lesssim \delta_{*}\langle t\rangle^{-\frac{1}{2}(\sigma-\sigma_0)}. 
\end{align*}
On the other hand, by decomposing  $u\cdot\nabla u=u^\ell\cdot\nabla u^\ell+u^h\cdot\nabla u^\ell+u^\ell\cdot\nabla u^h+u^h\cdot\nabla u^h$  for $\tau \geq 1$, we obtain that
\begin{align*} 
J_2\lesssim &\,   \int_1^t \langle t-\tau\rangle^{-\frac{1}{2}(\sigma-\sigma_0)} \Big(\|u\|_{{\dot B}_{2,1}^{\frac{d}{2}}}^\ell  \|u\|_{{\dot B}_{2,\infty}^{ \sigma_0+1}}^\ell+\|u\|_{{\dot B}_{2,1}^{\frac{d}{2}+1}}^h\|u\|_{{\dot B}_{2,\infty}^{ \sigma_0+1}}^\ell\nonumber\\
&+ \|u\|_{{\dot B}_{2,1}^{\frac{d}{2}}}^\ell \|u\|_{{\dot B}_{2,1}^{ \frac{d}{2}+1}}^h+\|u\|_{{\dot B}_{2,1}^{\frac{d}{2}}}^h \|u\|_{{\dot B}_{2,1}^{ \frac{d}{2}+1}}^h\Big){\rm d}\tau\nonumber\\
\lesssim&\, \mathcal{D}^2(t) \int_0^t  \langle t-\tau\rangle^{-\frac{1}{2}(\sigma-\sigma_0)} \langle\tau\rangle^{-\frac{1}{2}(\frac{d}{2}+1-\sigma_0)} {\rm d}\tau \lesssim  \mathcal{D}^2(t)\langle t\rangle^{-\frac{1}{2}(\sigma-\sigma_0)},
\end{align*}
for any $t\geq 2$. Here, we have utilized Lemma \ref{LA.8} with $\gamma_1 = \frac{1}{2}(\sigma - \sigma_0) \in \left[0, \frac{1}{2}\left(\frac{d}{2} + 1 - \sigma_0\right)\right)$ and $\gamma_2 = \frac{1}{2}\left(\frac{d}{2} + 1 - \sigma_0\right) > 1$.
By combining the estimates of $J_1$ and $J_2$, one gets   that for any $t\geq 2$,
\begin{align*}
\int_0^t \langle t-\tau\rangle^{-\frac{1}{2}(\sigma-\sigma_0)}\|u\cdot\nabla u\|_{\dot B_{2,\infty}^{\sigma_0} }^\ell{\rm d}\tau\lesssim     \mathcal{D}^2(t)\langle t\rangle^{-\frac{1}{2}(\sigma-\sigma_0)},
\end{align*}
which together with \eqref{NJKG4.31} yields
\begin{align}\label{NJKG4.32}
\int_0^t \langle t-\tau\rangle^{-\frac{1}{2}(\sigma-\sigma_0)}\|u\cdot\nabla u\|_{\dot B_{2,\infty}^{\sigma_0} }^\ell{\rm d}\tau\lesssim    ( \mathcal{D}^2(t)+\delta_*)\langle t\rangle^{-\frac{1}{2}(\sigma-\sigma_0)},    
\end{align}
for any $t>0$.

For the remaining terms, observe that no dissipative structure is associated with $a$. The coupling terms $f(a)\Delta u$ and $f(a)\nabla{\rm div}\,u$   lead to a loss of derivatives, in stark contrast to the situation in classical isentropic compressible Navier-Stokes equations, for example, see \cite{Danchin-IM-00}. 
 
\begin{rem}\label{NJKR4.1}
In fact, for the term $f(a)\Delta u$, by employing the energy method   as those presented in \cite[Section 5.2]{Danchin-2018} or \cite[Section 4.2]{LS-SIMA-2023}, we obtain  that
\begin{align*}
 &\int_1^t  \langle t-\tau\rangle^{-\frac{1}{2}(\sigma-\sigma_0)}\|f(a) \Delta u\|_{\dot B_{2,\infty}^{\sigma_0} }^\ell{\rm d}\tau\nonumber\\
& \lesssim\, \mathcal{X}(t)\int_0^t   \langle t-\tau\rangle^{-\frac{1}{2}(\sigma-\sigma_0)} \|u\|^\ell_{\dot B_{2,\infty}^{\sigma_0+2}} {\rm d}\tau+\mathcal{X}(t)\int_0^t   \langle t-\tau\rangle^{-\frac{1}{2}(\sigma-\sigma_0)} \|u\|^h_{\dot B_{2,\infty}^{\sigma_0+2}} {\rm d}\tau \nonumber\\
 &\lesssim\,\delta_{*}\int_0^t  \langle t-\tau\rangle^{-\frac{1}{2}(\sigma-\sigma_0)} \langle \tau\rangle^{-1}  {\rm d}\tau+\delta_*\langle t\rangle^{-\frac{1}{2}(\sigma-\sigma_0)},  
\end{align*}    
However, the integral term   {$\int^t_0\langle t-\tau\rangle^{-\frac{1}{2}(\sigma-\sigma_0)} \langle \tau\rangle^{-1}  {\rm d}\tau$} cannot be bounded by $\langle t\rangle^{-\frac{1}{2}(\sigma-\sigma_0)}$.
The argument for $f(a)\nabla{\rm div}\,u$ is similar.
\end{rem}
Therefore, these terms require careful handling in the analysis. Fortunately, we can handle them by adopting the similar approach outlined in \cite[Section 5]{LSZ-arXiv-2025}. By imposing additional conditions such as $a^\ell \in \dot B_{2,\infty}^{\sigma_0+1}$, we can effectively overcome the challenges arising from the absence of a dissipative structure for $a$.
It follows from \eqref{A4}$_1$ and {Lemma \ref{LA.7}} that
\begin{align*}
\|a\|_{\widetilde L_t^\infty(\dot B_{2,\infty}^{\sigma_0+1})}\lesssim&\,
\exp\Big\{C\| u\|_{L_t^1(\dot B_{2,1}^{\frac{d}{2}+1})} \Big\}
\Big( \|a_0\|_{\dot B_{2,\infty}^{\sigma_0+1}}+\|a{\rm div}\, u\|_{L_t^1(\dot B_{2,\infty}^{\sigma_0+1})} +\|{\rm div}\, u\|_{L_t^1(\dot B_{2,\infty}^{\sigma_0+1})}  \Big)  \nonumber\\
\lesssim&\,  {\|a_0\|_{\dot B_{2,\infty}^{\sigma_0+1}}}+\|a\|_{\widetilde L_t^\infty(\dot B_{2,\infty}^{\sigma_0+1})}\|u\|_{L_t^1(\dot B_{2,1}^{\frac{d}{2}+1})}+\|u\|_{L_t^1(\dot B_{2,\infty}^{\sigma_0+2})}^\ell+\|u\|_{L_t^1(\dot B_{2,\infty}^{\frac{d}{2}+1})}^h\nonumber\\
\lesssim&\, \delta_0\|a\|_{\widetilde L_t^\infty(\dot B_{2,\infty}^{\sigma_0+1})}+\delta_{*},
\end{align*}
which together with the smallness of $\delta_0$ yields
\begin{align}\label{NJKG4.33}
 \|a\|_{\widetilde L_t^\infty(\dot B_{2,\infty}^{\sigma_0+1})}\lesssim \delta_*.   
\end{align}

Next, we proceed to estimate the remaining nonlinear terms in \eqref{NJKG4.30}. For brevity, we focus on the non-trivial case where $t \geq 2$. Applying frequency decomposition, we have  
\begin{align*}
 &\int_1^t  \langle t-\tau\rangle^{-\frac{1}{2}(\sigma-\sigma_0)}\Big(\|f(a) \Delta u\|_{\dot B_{2,\infty}^{\sigma_0} }^\ell+{\|f(a) \Delta u\|_{\dot B_{2,\infty}^{\sigma_0} }^h}\Big){\rm d}\tau\nonumber\\
 \lesssim&\,  \int_0^t   \langle t-\tau\rangle^{-\frac{1}{2}(\sigma-\sigma_0)}\|a\|_{\dot B_{2,\infty}^{\sigma_0+1}} \Big(\|u\|^\ell_{\dot B_{2,1}^{\frac{d}{2}+1}}+\|u\|^h_{\dot B_{2,1}^{\frac{d}{2}+1}}\Big) {\rm d}\tau  \nonumber\\
 \lesssim&\, \delta_*\mathcal{D}(t)\int_0^t\langle t-\tau\rangle^{-\frac{1}{2}(\sigma-\sigma_0)} \langle\tau \rangle^{-\frac{1}{2}(\frac{d}{2}+1-\sigma_0)} {\rm d}\tau\nonumber\\
 \lesssim&\,\delta_*\mathcal{D}(t)\langle t\rangle^{-\frac{1}{2}(\sigma-\sigma_0)},
\end{align*}
where we have used \eqref{NJKG4.33} and Lemma \ref{LA.8}. Another case is evident and analogous to \eqref{NJKG4.31}. Consequently, we  infer that
\begin{align}\label{NJKG4.34}
 &\int_0^t  \langle t-\tau\rangle^{-\frac{1}{2}(\sigma-\sigma_0)}\Big(\|f(a) \Delta u\|_{\dot B_{2,\infty}^{\sigma_0} }^\ell+ {\|f(a) \Delta u\|_{\dot B_{2,\infty}^{\sigma_0} }^h}\Big){\rm d}\tau\lesssim \delta_*\big(1+\mathcal{D}(t)\big)\langle t\rangle^{-\frac{1}{2}(\sigma-\sigma_0)},
\end{align}
for any $t>0$.
Putting \eqref{NJKG4.32} and \eqref{NJKG4.34} into \eqref{NJKG4.30}, we further obtain that
\begin{align*}
\|u\|_{\dot B_{2,1}^\sigma}^\ell\lesssim&\,    (\delta_*+\delta_*\mathcal{D}(t)+\mathcal{D}^2(t))\langle t\rangle^{-\frac{1}{2}(\sigma-\sigma_0)},
\end{align*}
for any $t>0$, which  completes the proof of \eqref{NJKG4.28}.
\end{proof}
\begin{lem}\label{NJKL4.5}
Let $(a,u)$ be the global strong solution to the Cauchy problem \eqref{A4} given by Theorem \ref{T1.2}. Then, under the assumptions of Theorem \ref{T1.1},  it holds that for any $t>0$,
\begin{align}\label{NJKG4.35}
 \|\langle\tau\rangle^\alpha u\|_{\widetilde L_t^\infty(\dot B_{2,1}^{\frac{d}{2}-1})}^h+\| \tau ^\alpha u\|_{\widetilde L_t^\infty(\dot B_{2,1}^{\frac{d}{2}+1})}^h\lesssim  \delta_*+\delta_*\mathcal{D}(t),  
\end{align}
where $\delta_{*}$ and $\mathcal{D}(t)$ are defined via                 \eqref{TC.3} and \eqref{NJKG4.27}, respectively.  
\end{lem}
\begin{proof}
From \eqref{G3.12}, it holds that
\begin{align}\label{NJKG5.36}
\|\dot\Delta_k u\|_{L^2}\lesssim&\, e^{-t}\|\dot\Delta_k u_0\|_{L^2}+\int_0^t e^{-(t-\tau)}   \big(\|{\rm div}\, u\|_{L^\infty}\|\dot\Delta_k u\|_{L^2}+\|[u\cdot\nabla,\dot\Delta_k] u\|_{L^2}    \big){\rm d}\tau\nonumber\\
& +\int_0^te^{-(t-\tau)}\big(  \|\dot\Delta_k(f(a)\Delta u)\|_{L^2}+\|\dot\Delta_k(f(a)\nabla{\rm div}\, u)\|_{L^2}\big)  {\rm d}\tau,   
\end{align}
for any $k\geq -1$.
Applying the time-weighted method to \eqref{NJKG5.36} gives
\begin{align}\label{NJKG4.37}
\|\langle \tau\rangle^\alpha u \|_{\widetilde L_t^\infty(\dot B_{2,1}^{\frac{d}{2}-1})}^h\lesssim  \|u_0\|_{\dot B_{2,1}^{\frac{d}{2}-1}}^h+\sum_{k\geq -1}\sup_{\tau\in[0,t]}\langle\tau \rangle^\alpha \int_0^\tau e^{-(\tau-s)} 2^{k(\frac{d}{2}-1)}\sum_{j=1}^3K_j {\rm d}s,
\end{align}
where 
\begin{align*}
K_1=&\, \|{\rm div}\, u\|_{L^\infty}\|\dot\Delta_k u\|_{L^2},\\
K_2=&\, \|[u\cdot\nabla,\dot\Delta_k] u\|_{L^2}  ,\\
K_3=&\,  \|\dot\Delta_k(f(a)\Delta u)\|_{L^2}+\|\dot\Delta_k(f(a)\nabla{\rm div}\, u)\|_{L^2}.
\end{align*}
We primarily examine two cases: $t \leq 2$ and $t \geq 2$. For the case $t \leq 2$,  one has
\begin{align}\label{NJKG4.38}
&\sum_{k\geq -1}\sup_{\tau\in[0,t]}\langle\tau \rangle^\alpha \int_0^\tau e^{-(\tau-s)} 2^{k(\frac{d}{2}-1)} (K_1+K_2) {\rm d}s\nonumber\\
\lesssim&\, \int_0^t \|{\rm div}\,u\|_{L^\infty} \|u\|^h_{\dot B_{2,1}^{\frac{d}{2}-1}}  {\rm d}s+\int_0^t \| u\|_{\dot B_{2,1}^{\frac{d}{2}+1}} \|u\|_{\dot B_{2,1}^{\frac{d}{2}-1}}  {\rm d}s\nonumber\\ 
\lesssim&\,  \|u\|_{L_t^1(\dot B_{2,1}^{\frac{d}{2}+1})} \|u\|_{\widetilde L_t^\infty(\dot B_{2,1}^{\frac{d}{2}-1})}\lesssim \mathcal{X}^2(t).
\end{align}
As for the case $t\geq 2$, we divide the time interval $[0, t]$ into two parts: $[0, 1]$ and $[1, t]$. On the one hand, for $[0,1]$, the direct calculation yields 
\begin{align}\label{NJKG4.39}
&\sum_{k\geq -1}\sup_{\tau\in[2,t]}\langle\tau \rangle^\alpha \int_0^1 e^{-(\tau-s)} 2^{k(\frac{d}{2}-1)} (K_1+K_2) {\rm d}s  \lesssim   \mathcal{X}^2(1). 
\end{align} 
On the other hand,  for $[1,t]$, it holds that 
\begin{align}\label{NJKG4.40}
 &\sum_{k\geq -1}\sup_{\tau\in[2,t]}\langle\tau \rangle^\alpha \int_1^\tau e^{-(\tau-s)} 2^{k(\frac{d}{2}-1)} (K_1+K_2) {\rm d}s\nonumber\\
 &\lesssim\, \Big(\|\tau^\alpha u^\ell\|_{\widetilde L_t^\infty(\dot B_{2,1}^{\frac{d}{2}+1})}+\|\tau^\alpha u^h\|_{\widetilde L_t^\infty(\dot B_{2,1}^{\frac{d}{2}+1})}\Big) \| u\|_{\widetilde L_t^\infty(\dot B_{2,1}^{\frac{d}{2}-1})} \sup_{\tau\in [2,t]} \int_1^\tau e^{-(\tau-s)}s^{-\alpha} {\rm d}s\nonumber\\
 &\lesssim\, \mathcal{X}(t) \Big( \mathcal{D}(t)+ \|\tau^\alpha u^\ell\|_{\widetilde L_t^\infty(\dot B_{2,1}^{\frac{d}{2}+1})}         \Big).
\end{align}
Thanks to the definition of $\mathcal{D}(t)$, we get
\begin{align*}
 \|\tau^\alpha u^\ell\|_{\widetilde L_t^\infty(\dot B_{2,1}^{\frac{d}{2}+1})}\lesssim \mathcal{D}(t),
\end{align*}
which together with \eqref{NJKG4.40} gives
\begin{align}\label{NJKG4.41}
 &\sum_{k\geq -1}\sup_{\tau\in[2,t]}\langle\tau \rangle^\alpha \int_1^\tau e^{-(\tau-s)} 2^{k(\frac{d}{2}-1)} (K_1+K_2) {\rm d}s\lesssim \delta_*\mathcal{D}(t).   
\end{align}
Combining the estimates \eqref{NJKG4.38}, \eqref{NJKG4.39} and \eqref{NJKG4.41} up, we end up with
\begin{align}\label{NJKG4.42}
\sum_{k\geq -1}\sup_{\tau\in[0,t]}\langle\tau \rangle^\alpha \int_0^\tau e^{-(\tau-s)} 2^{k(\frac{d}{2}-1)} (K_1+K_2) {\rm d}s\lesssim \delta_*\mathcal{D}(t),    
\end{align}
for any $t>0$.
Similarly, it can be concluded that
\begin{align}\label{NJKG4.43}
\sum_{k\geq -1}\sup_{\tau\in[0,t]}\langle\tau \rangle^\alpha \int_0^\tau e^{-(\tau-s)} 2^{k(\frac{d}{2}-1)}{K_3} {\rm d}s\lesssim \delta_*\mathcal{D}(t),    
\end{align}
for any $t>0$.
Plugging the estimates \eqref{NJKG4.42} and \eqref{NJKG4.43} into \eqref{NJKG4.37}, we derive that
\begin{align}\label{NJKG4.44}
 \|\langle\tau\rangle^\alpha u\|_{\widetilde L_t^\infty(\dot B_{2,1}^{\frac{d}{2}-1})}^h\lesssim \delta_*+\delta_*\mathcal{D}(t) ,
\end{align}
for any $t>0$.
For the remaining high-regularity estimate of $u$, by taking $M = \alpha > 1$, we obtain that
\begin{align}\label{NJKG4.45}
\|\tau^\alpha u(t)\|_{\widetilde L_t^\infty( \dot B_{2,1}^{\frac{d}{2}+1})}^h
\lesssim&\,  \|\tau^{\alpha-1 } u\|_{\widetilde L_t^\infty( \dot B_{2,1}^{\frac{d}{2}-1})}^h   +\|\tau^{\alpha } u\cdot\nabla u\|_{\widetilde L_t^\infty( \dot B_{2,1}^{\frac{d}{2}-1})}^h
+\|\tau^{\alpha} f(a)\Delta u\|_{\widetilde L_t^\infty( \dot B_{2,1}^{\frac{d}{2}-1})}^h\nonumber\\
&+\|\tau^{\alpha} f(a)\nabla{\rm div}\,u\|_{\widetilde L_t^\infty( \dot B_{2,1}^{\frac{d}{2}-1})}^h\nonumber\\
\lesssim&\,  \|\langle\tau\rangle^\alpha u\|_{\widetilde L_t^\infty(\dot B_{2,1}^{\frac{d}{2}-1})}^h+  \|u\|_{\widetilde L_t^\infty( \dot B_{2,1}^{\frac{d}{2}-1})} \Big ( \|\tau^\alpha u\|_{\widetilde L_t^\infty( \dot B_{2,1}^{\frac{d}{2}+1})}^\ell+  \|\tau^\alpha u\|_{\widetilde L_t^\infty( \dot B_{2,1}^{\frac{d}{2}+1})}^h  \Big) \nonumber\\
&+ \|a\|_{\widetilde L_t^\infty(\dot B_{2,1}^\frac{d}{2})} \Big ( \|\tau^\alpha u\|_{\widetilde L_t^\infty( \dot B_{2,1}^{\frac{d}{2}+1})}^\ell+  \|\tau^\alpha u\|_{\widetilde L_t^\infty( \dot B_{2,1}^{\frac{d}{2}+1})}^h  \Big) \nonumber\\
\lesssim&\, \delta_*\mathcal{D}(t),    
\end{align}
for any $t>0$.
By combining \eqref{NJKG4.44} with
\eqref{NJKG4.45}, we obtain \eqref{NJKG4.35} directly.
\end{proof}

We now proceed to the proof of Theorem \ref{T1.1}.
\begin{proof}[Proof of Theorem \ref{T1.1}]
It follows from Lemmas \ref{NJKL4.4} and \ref{NJKL4.5} that
\begin{align*}
\mathcal{D}(t)\lesssim  \delta_*+\delta_*\mathcal{D}(t
)+\mathcal{D}^2(t),   
\end{align*}
for any $t>0$. The smallness of $\|u_0^\ell\|_{\dot B_{2,\infty}^{\sigma_0}}$ ensures that $\delta_*$ is sufficiently small. By exploiting the smallness of $\delta_*$, we consequently obtain $\mathcal{D}(t) \lesssim \delta_*$. Therefore, the estimates \eqref{TE.2} and \eqref{TE.4} hold, and the proof of Theorem \ref{T1.1} is completed.
\end{proof}

\subsection{The lower bound estimate for $u$}
In this subsection, we investigate the lower bound  of the decay estimates for $u$.
To begin with, we consider the following linearized system of \eqref{A4}$_2$:
\begin{equation}\label{G5.26}\left\{
\begin{aligned}            
& \partial_t u_{L} -\mu\Delta u_{L}-(\mu+\nu)\nabla {\rm div}\, u_{L}= 0, \\
&  u_{L}(0,x)=u_0(x).
\end{aligned}\right.
\end{equation}
Inspired by the approach developed in \cite{Danchin-IM-00}, we apply the orthogonal projectors $\mathbb{P}$ and $\mathbb{Q}$ to decompose 
the velocity   field $u_L$ into divergence-free and potential components, respectively. This allows the linearized system \eqref{G5.26} to be reformulated  into
\begin{equation}\label{G5.27} 
 \partial_t \mathbb{P}u_{L}- \mu\Delta\mathbb{P}u_{L}  = 0,\quad   \mathbb{P}u_{L}(0,x)=\mathbb{P}{u }_0(x), 
\end{equation}
and  
\begin{align}  \label{GN5.28} 
  \partial_t \mathbb{Q}u_{L}-(2\mu+\nu) \Delta\mathbb{Q}u_{L}  = 0, \quad  
 \quad \mathbb{Q}u_{L}(0,x)=\mathbb{Q}{u}_0(x), 
\end{align} 
where $\mathbb{P}u$ and $\mathbb{Q}u$ satisfy the ordinary heat equations.
Taking the Fourier transform to \eqref{G5.27}  and \eqref{GN5.28} with respect to $x$ gives
\begin{align*}
 \mathcal{F}[\mathbb{P}u_{L}](t,\xi)=e^{-\mu|\xi|^2 t} \mathcal{F}[\mathbb P{u}_0],\quad \mathcal{F}[\mathbb{Q}u_{L}](t,\xi)=e^{-(2\mu+\nu)|\xi|^2 t} \mathcal{F}[\mathbb Q{u}_0],
\end{align*}
  which implies that
\begin{align*}
 \widehat u_{L}(t,\xi)=&\,e^{-\mu|\xi|^2 t}\widehat u_0+\big(e^{-(2\mu+\nu)|\xi|^2 t}-e^{-\mu|\xi|^2 t}\big)\frac{\xi\xi^\top}{|\xi|^2} \widehat u_0\nonumber\\
 =:&\, \mathbb P(\xi) e^{-\mu|\xi|^2 t} \widehat u_0+\mathbb Q(\xi)e^{-(2\mu+\nu)|\xi|^2 t}\widehat u_0,
\end{align*}
where $\mathbb P(\xi):=\Big(1-\frac{\xi\xi^\top}{|\xi|^2}\Big)$ and
$\mathbb Q(\xi):=\frac{\xi\xi^\top}{|\xi|^2}$.
By virtue of   the orthogonality of $\mathbb P(\xi)$ and $\mathbb Q(\xi)$,   we further arrive at
\begin{align} \label{G5.28}
|\widehat u_{L}(t,\xi)|=&\,\Big(|\mathbb P(\xi)|^2e^{-2\mu|\xi|^2 t} |\widehat u_0|^2+|\mathbb Q(\xi)|^2e^{-2(2\mu+\nu)|\xi|^2 t} |\widehat u_0|^2 \Big)^{\frac{1}{2}}\nonumber\\
\gtrsim&\, e^{-\max\{\mu,2\mu+\nu\}}|\widehat{u}_0|.
 \end{align} 

Now, we state the decay estimate of the linearized system \eqref{G5.26} as follows.
\begin{lem}[Linear analysis]\label{L5.4}
  Let  $d\geq 2$.  Assume that $u_{L}$ is a  solution to
the   Cauchy problem \eqref{G5.26}, and the initial data $ u_0 $ satisfies
\begin{align*}
  u_0^\ell\in \dot {\mathcal B}_{2,\infty}^{\sigma_0}\quad \text{with}\quad \sigma_0\in\Big[-\frac{d}{2},\frac{d}{2}-1\Big),
\end{align*}  
then for all $t \geq 1$, there exists two universal constants $c_5>0$ and  $C_5 > 0$ such that  
\begin{align} \label{G5.29}
c_5(1+t)^{-\frac{1}{2}(\sigma-\sigma_0)}\leq \|u_{L}(t)\|_{\dot B_{2,1}^\sigma}\leq C_5  (1+t)^{-\frac{1}{2}(\sigma-\sigma_0)},    
\end{align}
for all $\sigma\in(\sigma_0,\frac{d}{2}+1]$, where $\dot{\mathcal{B}}_{2,\infty}^{\sigma_0}$ is defined by \eqref{G1.12}.
\end{lem}
\begin{proof}
The upper bound of the decay estimate for $u_L$ follows directly from Theorem \ref{T1.4}. For brevity, we omit the details here. Next, we focus exclusively on establishing the lower bound estimate for $u_L$ in the system \eqref{G5.26}, where the proof is inspired by \cite{Bl-2016-SIMA} and \cite[Section 3]{BSXZ-Adv-2024}. Without loss of generality, for $\{k_j\}_{j\in\mathbb{N}} \subset \mathbb{Z}$, we assume that $j = 1, 2, \dots$ correspond to the indices less than $[\log_2 \zeta]$. By leveraging the Fourier-Plancherel theorem and \eqref{G5.29}, we have
\begin{align*}
\|u_L(t)\|_{\dot B_{2,1}^\sigma}\geq  \|u_L^\ell(t)\|_{\dot B_{2,1}^\sigma}  \geq&\,  \sum_{k\leq [\rm log_2\zeta]} 2^{\sigma k} \|\dot\Delta_k u_L(t)\|_{L^2}   \nonumber\\
\geq&\, \sum_{k\leq [\rm log_2\zeta]}2^{\sigma k} \|\phi (2^{-k}\cdot)\widehat u_L(t)\|_{L^2{}}\nonumber\\
\geq&\, \sum_{k\leq [\rm log_2\zeta]}e^{-\frac{64\max\{\mu,2\mu+\nu\} 2^{2k}t}{9}}2^{\sigma k} \|\dot\Delta_k u_0\|_{L^2},
\end{align*}
for any fixed $t\geq 1$. From \eqref{G1.12}, we can find a maximal integer $k_{j_0}$ satisfying $k_{j_0} \leq -\frac{1}{2} \log_2(1+t)$. We claim that $k_{j_0} > -M_0 - \frac{1}{2} \log_2(1+t)$. Otherwise, if there exists another integer $k_{j_0+1}$ such that $k_{j_0+1} \leq k_{j_0} + M_0 \leq -\frac{1}{2} \log_2(1+t)$, this would contradict the maximality of $k_{j_0}$. It follows from \eqref{G1.12} and $2^{k_{j_0}}\backsim (1+t)^{-\frac{1}{2}}$ that
\begin{align*}
\|u_L(t)\|_{\dot B_{2,1}^\sigma} \gtrsim  \|u_L^\ell(t)\|_{\dot B_{2,1}^\sigma} 
\gtrsim&\, \sum_{k\leq [\rm log_2\zeta]}e^{-\frac{64\max\{\mu,2\mu+\nu\} 2^{2k}t}{9}}2^{\sigma k} \|\dot\Delta_k u_0\|_{L^2}\nonumber\\
\gtrsim&\, e^{-\frac{64\max\{\mu,2\mu+\nu\} 2^{2 k_{j_0}} t  } {9}} 2^{(\sigma-\sigma_0)k_{j_0}} 2^{\sigma_0 k_{j_0}} \|\dot\Delta_k u_0\|_{L^2}\nonumber\\
\gtrsim&\,  c_02^{(\sigma-\sigma_0)k_{j_0}}\nonumber\\
\gtrsim&\, c_0(1+t)^{-\frac{1}{2}(\sigma-\sigma_0)},
\end{align*}
for all $t\geq 1$. Therefore, we  prove \eqref{G5.27}, thereby completing the proof of Lemma \ref{L5.4}.
\end{proof}

To study the nonlinear component of \eqref{A4}$_2$, we define $\omega = u - u_L$ and then consider the following nonlinear Cauchy problem associated with $\omega$:
 \begin{equation}\label{G5.30}\left\{
\begin{aligned}            
& \partial_t \omega-\mu\Delta \omega-(\mu+\nu)\nabla {\rm div}\, \omega = F, \\
& \omega(0,x)=0,
\end{aligned}\right.
\end{equation}
where $F:= -u\cdot\nabla u+ f(a)\mu\Delta u+f(a)(\mu+\nu)\nabla{\rm div}\, u$.

\begin{lem}[Nonlinear analysis]\label{L4.7}
Let $d\geq 2$. Assume that $\omega$ is
a solution to the nonlinear Cauchy problem \eqref{G5.30},
and the initial data $(a_0, u_0)$  defined in \eqref{A4}  satisfies
\begin{align*}
a^\ell_0\in \dot B_{2,\infty}^{\sigma_0+1},\quad u_0^\ell\in \dot {\mathcal B}_{2,\infty}^{\sigma_0},\quad \|u_0^\ell\|_{\dot B_{2,\infty}^{\sigma_0}}\leq\eps_3,\quad \text{with}\quad \sigma_0\in\Big[-\frac{d}{2},\frac{d}{2}-1\Big),
\end{align*}
where $\eps_3$  is sufficiently small positive constant, 
then  it holds that for all $t > 0$,
\begin{align*} 
\|\omega(t) \|_{\dot B_{2,1}^\sigma}   \lesssim \delta_{*}^2  (1+t)^{-\frac{\sigma-\sigma_0}{2}},
\end{align*}
for all $\sigma\in(\sigma_0,\frac{d}{2}+1]$, where $\dot{\mathcal{B}}_{2,\infty}^{\sigma_0}$ is defined via \eqref{G1.12}.
\end{lem}
\begin{proof}
Similar to \eqref{G3.6}, we get
\begin{align} \label{G5.31}
&\frac{{\rm d}}{{\rm d}t}\|\dot\Delta_k \omega\|_{L^2}^2+2^{2k}\| \dot\Delta_k \omega\|_{L^2}^2 \lesssim \|\dot\Delta_k\omega\|\|\dot\Delta_kF\|_{L^2},
\end{align}
for any $k\leq 0$.
Applying Gronwall's inequality to \eqref{G5.31}, together with the initial condition $\omega(0,x)=0$, yields
\begin{align*} 
\|\dot\Delta_k\omega\|_{L^2}\lesssim e^{-2^{2k}t}   \|\dot \Delta_k\omega_0\|_{L^2}+\int_0^t e^{-2^{2k}(t-\tau)}\|\dot\Delta_k F\|_{L^2}{\rm d}\tau\lesssim \int_0^t e^{-2^{2k}(t-\tau)}\|\dot\Delta_k F\|_{L^2}{\rm d}\tau , 
\end{align*}
which implies that for any $\sigma\in(\sigma_0,\frac{d}{2}-1]$,
\begin{align*}
 \|\omega\|_{\dot B_{2,1}^{\sigma}}^\ell \lesssim&\, \int_0^t \langle t-\tau\rangle^{-\frac{1}{2}(\sigma-\sigma_0)} \|F\|_{\dot B_{2,\infty}^{\sigma_0}}^\ell {\rm d}\tau \nonumber\\
 \lesssim&\,  \int_0^t \langle t-\tau\rangle^{-\frac{1}{2}(\sigma-\sigma_0)} \Big(\|u\cdot\nabla u\|_{\dot B_{2,\infty}^{\sigma_0}}^\ell+\|f(a)\Delta u\|_{\dot B_{2,\infty}^{\sigma_0}}^\ell+\|f(a)\nabla {\rm div}\, u\|_{\dot B_{2,\infty}^{\sigma_0}}^\ell     \Big) {\rm d}\tau.
\end{align*}
Taking the same argument as that in Lemma \ref{NJKL4.4}, we   deduce that
\begin{align*}
\|\omega\|_{\dot B_{2,1}^{\sigma}}^\ell \lesssim ( \delta_*\mathcal{D}(t)+\mathcal{D}^2(t))\langle t\rangle^{-\frac{1}{2}(\sigma-\sigma_0)}\lesssim \delta_*^2(1+t)^{-\frac{1}{2}(\sigma-\sigma_0)},    
\end{align*}
for any $t\geq 1$.
\end{proof}
With the Lemmas \ref{L5.4} and \ref{L4.7} in hand, we are now in a position to prove Theorem \ref{T1.5}.
\begin{proof}[Proof of Theorem \ref{T1.5}]
By  Duhamel’s principle, we have that for any $\sigma\in (\sigma_0,\frac{d}{2}+1]$,
\begin{align*}
\|u(t)\|_{\dot B_{2,1}^\sigma} \geq \|u^\ell(t)\|_{\dot B_{2,1}^\sigma}\geq&\, \|u^\ell_L(t)\|_{\dot B^\sigma_{2,1}} -\|\omega(t)\|_{\dot B_{2,1}^\sigma}^\ell\nonumber\\
\geq&\, c_0 (1+t)^{-\frac{1}{2}(\sigma-\sigma_0)}-\delta_*^2(1+t)^{-\frac{1}{2}(\sigma-\sigma_0)}\nonumber\\
\geq&\, c_6(1+t)^{-\frac{1}{2}(\sigma-\sigma_0)},
\end{align*}
for all $t \geq 1$ and some constant $c_6 > 0$. On the other hand, it is straightforward to verify that, from Theorem \ref{T1.1}, for any $\sigma \in (\sigma_0, \frac{d}{2}+1]$, 
\begin{align*}
\|u(t)\|_{\dot B_{2,1}^\sigma} \lesssim \delta_* (1+t)^{-\frac{1}{2}(\sigma - \sigma_0)},
\end{align*}
for all $t \geq 1$. These inequalities complete the proof of Theorem \ref{T1.5}.
\end{proof}

\appendix 
\section{Analytic tools}
 This section is devoted to presenting  several fundamental properties of Besov spaces and product estimates, which have been frequently used in previous sections. It is worth noting that these properties also apply to Chemin-Lerner type spaces, provided that the time exponent satisfies H\"older's inequality with respect to the time variable. We first introduce the classical Bernstein inequalities.
\begin{lem}\label{LA.1}{\rm(\!\!\cite{BCD-Book-2011})}
Let $0<r<R$, $1\leq p\leq q\leq \infty$ and $m\in \mathbb{N}$. Define the ball $\mathcal{B}=\{\xi\in\mathbb{R}^{3}~| ~|\xi|\leq  R\}$ and the annulus $\mathcal{C}=\{\xi\in\mathbb{R}^{3}~|~ \lambda r\leq |\xi|\leq \lambda R\}$ . For any $f\in L^p$ and $\lambda>0$, it holds that
\begin{equation*} \left\{
\begin{aligned}
{\rm{supp}}\, \mathcal{F}(f) \subset& \lambda \mathcal{B} \Rightarrow \|D^{m}f\|_{L^q}\lesssim\lambda^{m+d(\frac{1}{p}-\frac{1}q{})}\|f\|_{L^p}, \nonumber\\
{\rm{supp}}\, \mathcal{F}(f) \subset& \lambda \mathcal{C} \Rightarrow \lambda^{m}\|f\|_{L^{p}}\lesssim\|D^{m}f\|_{L^{p}}\lesssim \lambda^{m}\|f\|_{L^{p}}.
\end{aligned}\right.
\end{equation*}
 \end{lem}
With the help of the above Bernstein inequalities in Lemma \ref{LA.1}, we can deduce the following properties of Besov spaces.
\begin{lem}\label{LA.2}{\rm(\!\!\cite[Chapter 2]{BCD-Book-2011})}
The following properties hold  :
\begin{itemize}
\item{} For $s\in\mathbb{R}$, $1\leq p_{1}\leq p_{2}\leq \infty$ and $1\leq r_{1}\leq r_{2}\leq \infty$, it holds
\begin{equation}\notag
\begin{aligned}
\dot{B}^{s}_{p_{1},r_{1}}\hookrightarrow \dot{B}^{s-d (\frac{1}{p_1} -\frac{1}{p_2})}_{p_{2},r_{2}}.
\end{aligned}
\end{equation}
\item{} For $1\leq p\leq q\leq\infty$, we have the following chain of continuous embedding  :
\begin{equation}\nonumber
\begin{aligned}
  \dot{B}^{0}_{p,1}\hookrightarrow L^{p}\hookrightarrow \dot{B}^{0}_{p,\infty}\hookrightarrow \dot B_{q,\infty}^{\varrho},\quad \varrho=-d\Big(\frac{1}{p}-\frac{1}{q}\Big).
\end{aligned}
\end{equation}
\item{} If $p<\infty$, then $\dot{B}^{\frac{d}{p}}_{p,1}$ is continuously embedded in the set of continuous functions decaying to 0 at infinity.

\item{}  The following real interpolation property is satisfied for $1\leq p\leq \infty$, $s_1<s_2$, and $\theta\in (0,1)$ :
\begin{align}\label{A.1}
\|f\|_{\dot B_{p,1}^{\theta s_1+(1-\theta)s_2}}\lesssim \frac{1}{\theta(1-\theta)(s_2-s_1)}\|f\|_{\dot B_{p,\infty}^{s_1}}^{\theta} \|f\|_{\dot B_{p,1}^{s_2}}^{1-\theta}  .  
\end{align}

\item{}  For any $\varepsilon>0$, it holds that
\begin{equation}\nonumber
\begin{aligned}
H^{s+\varepsilon}\hookrightarrow \dot{B}^{s}_{2,1}\hookrightarrow \dot{H}^{s}.
\end{aligned}
\end{equation}
\item{}
Let $\Lambda^{\sigma}$ be defined by $\Lambda^{\sigma}=(-\Delta )^{\frac{\sigma}{2}}f:=\mathcal{F}^{-1}\big{(} |\xi|^{\sigma}\mathcal{F}(f) \big{)}$ for $\sigma\in \mathbb{R}$ and $f\in{\mathcal S}^{'}_{h}(\mathbb{R}^3)$, then $\Lambda^{\sigma}$ is an isomorphism from $\dot{B}^{s}_{p,r}$ to $\dot{B}^{s-\sigma}_{p,r}$.
\item{} Let $1\leq p_{1},p_{2},r_{1},r_{2}\leq \infty$, $s_{1}\in\mathbb{R}$ and $s_{2}\in\mathbb{R}$ satisfy
\begin{align*}
 s_{2}<\frac{d}{p_2}\quad\text{\text{or}}\quad s_{2}=\frac{d}{p_2}~\text{and}~r_{2}=1.
\end{align*}
    The space $\dot{B}^{s_{1}}_{p_{1},r_{1}}\cap \dot{B}^{s_{2}}_{p_{2},r_{2}}$ endowed with the norm $\|\cdot \|_{\dot{B}^{s_{1}}_{p_{1},r_{1}}}+\|\cdot\|_{\dot{B}^{s_{2}}_{p_{2},r_{2}}}$ is a Banach space and has the weak compact and Fatou properties$:$ If $f_{n}$ is a uniformly bounded sequence of $\dot{B}^{s_{1}}_{p_{1},r_{1}}\cap \dot{B}^{s_{2}}_{p_{2},r_{2}}$, then an element $f$ of $\dot{B}^{s_{1}}_{p_{1},r_{1}}\!\cap \dot{B}^{s_{2}}_{p_{2},r_{2}}$ and a subsequence $f_{n_{k}}$ exist such that $f_{n_{k}}\rightarrow f$ in $\mathcal{S}'$ and
    \begin{align*}
    \begin{aligned}
    \|f\|_{\dot{B}^{s_{1}}_{p_{1},r_{1}}\cap \dot{B}^{s_{2}}_{p_{2},r_{2}}}\lesssim \liminf_{n_{k}\rightarrow \infty} \|f_{n_{k}}\|_{\dot{B}^{s_{1}}_{p_{1},r_{1}}\cap \dot{B}^{s_{2}}_{p_{2},r_{2}}}.
    \end{aligned}
    \end{align*}
\end{itemize}
\end{lem}

To control the nonlinear terms, we require the following Morse-type product estimates in Besov spaces:
\begin{lem} \label{LA.3}{\rm(\!\!\cite[Chapter 2]{BCD-Book-2011})}
The following statements hold  {\rm:}
\begin{itemize}
\item{} Let $s>0$, $1\leq p,r\leq \infty$.  Then $\dot B^{s}_{p,r}\cap L^\infty$ is an algebra and
 \begin{align}\label{A.2}
\|fg\|_{\dot B_{p,r}^s}\lesssim \|f\|_{L^\infty}\|g\|_{\dot B_{p,r}^s}+\|g\|_{L^\infty}\|f\|_{\dot B_{p,r}^s}.    
 \end{align}

\item{} Let $s_1,s_2>0$ and $p$ satisfy $2\leq p\leq \infty$, $s_1\leq \frac{d}{p}$, $s_2\leq \frac{d}{p}$, and $s_1+s_2>0$.  Then it holds 
 \begin{align}\label{A.3}
\|fg\|_{\dot B_{p,1}^{s_1+s_2-\frac{d}{p}}}\lesssim \|f\|_{ \dot B_{p,1}^{s_1 }}\|g\|_{\dot B_{p,1}^{s_2}} .    
 \end{align}

\item{} Let $s_1,s_2>0$ and $p$ satisfy $2\leq p\leq \infty$, $s_1\leq \frac{d}{p}$, $s_2<\frac{d}{p}$, and $s_1+s_2\geq0$.  Then it holds 
 \begin{align}\label{A.4}
\|fg\|_{\dot B_{p,\infty}^{s_1+s_2-\frac{d}{p}}}\lesssim \|f\|_{ \dot B_{p,1}^{s_1 }}\|g\|_{\dot B_{p,\infty}^{s_2}} .    
 \end{align}
    \end{itemize}
\end{lem}
We now present the following lemma concerning the continuity of composite functions:

\begin{lem}\label{LA.4}{\rm(\!\!\cite[Chapter 2]{BCD-Book-2011})}
Let $F:I\rightarrow \mathbb{R}$ be a smooth function such that $F(0)=0$. Then, for any $1\leq p\leq \infty$, $s>0$, and $f\in \dot{B}_{2,1}^s \cap L^\infty$, it holds that $F(f) \in \dot{B}_{2,1}^s \cap L^\infty$, and
\begin{align}\label{A.5}
\|F(f)\|_{\dot{B}^{s}_{2,1}} \leq C_f \|f\|_{\dot{B}_{2,1}^s},
\end{align}
where the constant $C_f > 0$ depends only on $\|f\|_{L^\infty}$, $F'$, $s$, and the spatial dimension $d$.
\end{lem}

To control the nonlinearities in high frequencies, we provide the following commutator estimates.
\begin{lem}\label{LA.5}
Let $1\leq p\leq \infty$ and $-\frac{d}{p}<s\leq 1+\frac{d}{p}$. Then it holds that
\begin{align}\label{A.6}
\sum_{k\in\mathbb Z} 2^{ks}\|[g\cdot\nabla,\dot\Delta_k]f\|_{L^p}\lesssim \|g\|_{\dot B^{\frac{d}{p}+1}_{p,1}}\|f\|_{\dot B_{p,1}^s},    
\end{align}
where the commutator $[A,B]:=AB-BA$.
\end{lem}

To achieve the optimal time decay rate of the solution, we present the following inequality:
\begin{lem}\label{LA.8} {\rm(\!\!\cite[Section 5]{Danchin-2018})}
Let $0<\gamma_1\leq\gamma_2$. If in addition $\gamma_2>1$, then it holds that
\begin{align}\label{A.10}
\int_0^t \langle t-\tau\rangle^{-\gamma_1} \langle \tau\rangle^{-\gamma_2} {\rm d}\tau\lesssim \langle  t\rangle^{-\gamma_1}.    
\end{align}
\end{lem}

Below, we consider the optimal regularity estimates for the Lamé system
\begin{align} \label{A.7}
\left\{\begin{aligned}
&\partial_t g-\mu\Delta  g-(\mu+\nu)\nabla{\rm div}\,g=f, \quad\, x\in \mathbb{R}^d,\quad t>0,  \\
& g(x,0)=g_0(x),\quad \,\,\,\quad\quad\quad\quad\quad\quad\quad  x\in \mathbb{R}^d.
 \end{aligned}
 \right.
\end{align}
\begin{lem}\label{LA.6}{\rm(\!\!\cite[Chapter 3]{BCD-Book-2011})}
Let $ T > 0 $, $ \mu > 0 $, $ 2\mu + \nu > 0 $, $ s \in \mathbb{R} $, $ 1 \leq p, r \leq \infty $, and $ 1 \leq \varrho_2 \leq \varrho_1 \leq \infty $. Suppose that $ g_0 \in \dot{B}^{s}_{p,r} $ and $ f \in \widetilde{L}^{\varrho_2}(0,T;\dot{B}_{p,r}^{s-2+\frac{2}{\varrho_2}}) $. Then there exists a solution $g $ to \eqref{A.7} satisfying
\begin{align*}
\min\{\mu,2\mu+\nu\}^{\frac{1}{\varrho_1}}\|g\|_{\widetilde{L}^{\varrho_1}_T(\dot B^{s+\frac{2}{\varrho_1}}_{p,r})}\lesssim \|g_0\|_{\dot B_{p,r}^s}+\min\{\mu,2\mu+\nu\}^{\frac{1}{\varrho_2}-1}\|f\|_{L^{\varrho_2}_T(\dot B_{p,r}^{s-2+\frac{2}{\varrho_2}})}.
\end{align*}
\end{lem}

Finally,   we consider the estimate of the transport equation
\begin{align} \label{A.8}
\left\{\begin{aligned}
&\partial_t h+g\cdot\nabla h=f, \quad\, x\in \mathbb{R}^d,\quad t>0,  \\
& h(x,0)=h_0(x),\quad \,\,\,\,\,  x\in \mathbb{R}^d.
 \end{aligned}
 \right.
\end{align}

\begin{lem}\label{LA.7}{\rm(\!\!\cite[Chapter 3]{BCD-Book-2011})}
Let $ T > 0 $, $ -\frac{d}{2} < s \leq \frac{d}{2} + 1 $, $ 1 \leq r \leq \infty $, $ h_0 \in \dot{B}_{2,r}^s $, $ g \in L^1(0,T; \dot{B}_{2,1}^{\frac{d}{2}+1}) $, and $ f \in L^1(0,T; \dot{B}_{2,r}^s) $. Then there exists a constant $ C > 0 $, independent of $ T $ and $ h_0 $, such that the solution $ h $ to \eqref{A.8} satisfies
\begin{align}\label{A.9}
\|h\|_{\widetilde L_T^\infty(\dot B_{2,r}^s)}\leq {\rm exp}\Big\{C\|\nabla g\|_{L_t^1(\dot B_{2,1}^{\frac{d}{2}})}\Big\}\bigg( \|h_0\|_{\dot B_{2,r}^s}+\int_0^t \|f\| _{\dot B_{2,r}^s}{\rm d}s\bigg) .   
\end{align}
Moreover, if $r<\infty$, then the solution $h$ belongs to $\mathcal{C}([0,T];\dot B_{2,r}^s)$.
\end{lem}
 
\bigskip 
{\bf Acknowledgements:}
Li and Ni are supported by NSFC (Grant No. 12331007 ).  
And Li is also supported by the ``333 Project" of Jiangsu Province.
Zhang  is supported by NSFC (Grant Nos. 12101305 and 12471215) and Taishan Scholars Program (tsqn202507101).

\vspace{2mm}

\textbf{Conflict of interest.} The authors do not have any possible conflicts of interest.

\vspace{2mm}

\textbf{Data availability statement.}
 Data sharing is not applicable to this article as no data sets were generated or analyzed during the current study.

\bibliographystyle{plain}

\end{document}